\documentclass{elsarticle}
\usepackage{graphicx}
\usepackage{amsfonts}
\usepackage{amsmath,mathtools}
\usepackage{verbatim}
\usepackage{amsthm}
\usepackage{tabularx,array}

\makeatletter
\def\ps@pprintTitle{%
 \let\@oddhead\@empty
 \let\@evenhead\@empty
 \def\@oddfoot{\centerline{\thepage}}%
 \let\@evenfoot\@oddfoot}
\makeatother
\providecommand{\OO}[1]{\mathop{\mathrm{O}}\bigl(#1\bigr)}
\providecommand{\oo}[1]{\mathop{\mathrm{o}}\bigl(#1\bigr)}

\newtheorem{prop}{Proposition}
\newtheorem{theorem}[prop]{Theorem}

\newtheorem{corollary}[prop]{Corollary}
\newtheorem{defin}[prop]{Definition}

\newtheoremstyle{rem} % name 
{3pt}% Space above 
{3pt}% Space below 
{}% Body font 
{}% Indent amount
{\bfseries}% Theorem head font 
{.}% Punctuation after theorem head 
{.5em}% Space after theorem head
{}% Theorem head spec (can be left empty, meaning `normal')
\theoremstyle{rem} 
\newtheorem{remark}[prop]{Remark}

\begin{document}

\begin{frontmatter}

\title{A probabilistic proof of cutoff in the Metropolis algorithm for the Erd\H{o}s-R\'{e}nyi random graph}
\tnotetext[t1]{This work is based on the PhD thesis \cite{Barta2012} of the author at the Department of Statistics, University of Chicago, 2012.}
\tnotetext[t2]{Current address: Department of Statistics, George Washington University, Washington, DC. Email address: barta@gwu.edu (Winfried Barta).}

\author{Winfried Barta}
%\ead{barta@gwu.edu}
\address{University of Chicago}

\begin{abstract}
We study mixing of the Metropolis algorithm for a distribution on the hypercube that corresponds to the Erd\H{o}s-R\'enyi random graph with edge probability $p$. This Markov chain has cutoff at $\max\{p,1-p\} n \log n$ with window size $n$, a result proved by Diaconis and Ram (2000) using Fourier analysis. Here we give an alternative proof that relies on coupling and a projection to a two-dimensional Markov chain. This is done in the hope that probabilistic techniques will be easier to generalize to less symmetric distributions. We also describe a close relationship between the Metropolis and Gibbs samplers for this model. 
Our proof extends to the case where the edge probabilities vary with $n$. In that case, we also show that a natural coordinate wise coupling is sharp if and only if 
%$\theta_n = \OO{1/n}$.
% Note that for $p:=\theta/(1+\theta)$ we have $\theta = \OO{1/n}$ if and only if $p = \OO{1/n}$.
the edge probabilities are of order $1/n$.
\end{abstract}


\end{frontmatter}

\section{Introduction}
We are interested in analyzing convergence rates of the random walk Metropolis algorithm for various distributions $\pi$ on the hypercube $\mathcal{X} := \{0,1\}^n$. This Markov chain on $\mathcal{X}$ moves as follows: Given that we are at a state $x \in \mathcal{X}$, we chose one of the $n$ neighbors of $x$ uniformly at random, say $y$, and propose to move from $x$ to $y$. Here $x$ is called a neighbor of $y$, denoted by $x \sim y$,  whenever $x$ and $y$ differ in exactly one coordinate. This proposal gets accepted, i.e. we move to $y$, with probability min$\big(1,\pi(y) / \pi(x)\big)$. If it gets rejected, we stay at $x$. This transition rule ensures that we have detailed balance
$$
\pi(x) \, P(x,y) = \pi(y) \, P(y,x) \qquad\mbox{for all $x,y \in \mathcal{X}$} .
$$
If $\pi$ is positive on the entire state space, then the chain is irreducible and its unique stationary distribution is $\pi$.
Often it will be convenient to make the chain \emph{lazy}. This corresponds to flipping a fair coin independently at each step. If it comes up heads, we stay where we are; if it comes up tails, we move according to the rule specified above.

We are interested in distributions $\pi$ that are \emph{unimodal} and \emph{radially symmetric} with respect to their mode.
By this we mean that there exists a state, $z$ say, that has highest mass under $\pi$;
any two states with the same distance to $z$ have the same $\pi$-mass; and this $\pi$-mass decreases with distance to $z$.
That is, we have
$$
\pi(x) = \pi(y) \qquad\mbox{whenever $d(x,z)=d(y,z)$} ,
$$
and
$$
\pi(x) \leq \pi(y) \qquad\mbox{whenever $d(x,z) > d(y,z)$}  .
$$
Here, $d(x,z) := \sum_{i=1}^n |x_i - z_i|$ is the graph distance (Hamming distance) of $x$ and $z$, i.e the number of coordinates where $x$ and $z$ differ. 

By relabeling the states, we may (and will) assume that the mode $z$ of $\pi$ is at $\underline{0}=(0,0,...,0)$. This ensures that $\pi$ is constant on level sets $L(k) := \{x \in \mathcal{X}: S(x)=k\}$ where $k \in \{0,1,...,n\}$ and $S(x):=\sum_{i=1}^n x_i$ is the number of ones in $x$. 
Hence, in the Metropolis algorithm we will always accept downward moves $x \rightarrow y$ where $S(x)>S(y)$, and we will accept upward moves $x \rightarrow y$ where $S(x)<S(y)$ with probability
$\theta_{S(x)}$, where we write $\theta_k := \pi(v)/\pi(w)$, where $v,w$ are some (any) states such that $S(v)=k+1$ and $S(w)=k$. 

For concreteness, the transition kernel for the lazy random walk Metropolis Hastings 
algorithm is
\begin{equation} \label{transition_kernel_lRW-MH}
P(x,y) =  \left\{ \begin{array}{l@{\quad:\quad}l}
\frac{1}{2n} & x \sim y, S(x) > S(y)  , \\
\frac{1}{2n} \theta_{S(x)} & x \sim y, S(x) < S(y) , \\
\frac{1}{2} + \frac{n-S(x)}{2n} (1-\theta_{S(x)}) & x = y , \\
0 & \mbox{otherwise}  .\\ 
\end{array} \right.
\end{equation}
Note that a consequence of radial symmetry is that the projection $S(X_t)$ of the Markov chain $(X_t)$ is also Markov, since for all $x,y$ we then have
$$
P(x,[y]) = P(x',[y]) \qquad\mbox{for all } x' \sim_S x .
$$
Here we write $x \sim_S y$ whenever $S(x)=S(y)$ and $[y] := \{z: S(z)=S(y)\}$ for $y \in \mathcal{X}$ are the equivalence classes of the relation $\sim_S$. 
See \cite[Lemma 2.5 on page 25]{LevinPeresWilmer2009}.

We measure distance to stationarity by total variation:
$$
||P^t(x,\cdot) - \pi||_{TV} := \max_{A \subset \mathcal{X}} \; \left(P^t(x,A)-\pi(A)\right)  ,
$$
and we are interested in this distance from the worst starting point:
$$
d(t) := \max_{x \in \mathcal{X}} ||P^t(x,\cdot) - \pi|| .
$$
The mixing time for a parameter $\varepsilon \in (0,1)$ is defined as
$$
t_{\rm mix}(\varepsilon) := \min\{t \geq 0: d(t) \leq \varepsilon\}
$$
and we write $t_{\rm mix}$ for $t_{\rm mix}(1/4)$. We are interested in the behavior of $t_{\rm mix}$ as the dimension $n$ of the hypercube goes to infinity. 

\vspace{0.5em}

An interesting phenomenon is that for some chains the mixing time  $t_{\rm mix}(\varepsilon)$ doesn't depend on the parameter $\varepsilon$ (asymptotically, as $n$ goes to infinity). We say a sequence $\left(\mathbf{X}^{(n)}\right)_{n \in \mathbb{N}}$ of Markov chains $\mathbf{X}^{(n)} = \big(X_t^{(n)}\big)_{t=0,1,...}$ on $\{0,1\}^n$ has a \emph{cutoff} 
(at $t^{(n)}_{\rm mix})$ 
if, for all $\varepsilon \in (0,1)$, 
\begin{equation}
\lim_{n \rightarrow \infty} \frac{t^{(n)}_{\rm mix}(\varepsilon)}{t^{(n)}_{\rm mix}(1-\varepsilon)} = 1  .
\end{equation}
Here, $t_{\rm mix}^{(n)}(\varepsilon)$ denotes the $\varepsilon$-mixing time of the $n^{th}$ chain $(X_t^{(n)})_{t=0,1,...}$. This is equivalent to 
\begin{equation}
\label{lemma18.1}
\lim_{n \rightarrow \infty} d_n\big(c \, t_{\rm mix}^{(n)}\big) = \left\{ \begin{array}{l@{\quad:\quad}l} 
1 & \mbox{if } \; c < 1  , \\
0 & \mbox{if } \; c > 1  . 
\end{array}
\right.
\end{equation}
So the function $d_n(\cdot)$, the total variation distance to stationarity from the worst starting point for the $n^{th}$ chain, approaches a step function as $n$ goes to infinity (if we rescale time by $t_{\rm mix}^{(n)}$). For a proof of this equivalence see \cite[Lemma 18.1 on page 247]{LevinPeresWilmer2009}, from where we also borrow the notation. For an overview of the cutoff phenomenon, see \cite{Diaconis1996}.

Sometimes it is possible to analyze more precisely what happens for $c=1$ in (\ref{lemma18.1}). We say a sequence of Markov chains has a cutoff with \emph{window size} $(w_n)$, if $w_n \in o(t_{\rm mix}^{(n)})$ and
\begin{eqnarray}
\lim_{\alpha \rightarrow \infty} \, \liminf_{n \rightarrow \infty} \; d_n\big(t_{\rm mix}^{(n)} - \alpha w_n\big) & = & 1  , \\
\lim_{\alpha \rightarrow \infty} \, \limsup_{n \rightarrow \infty} \; d_n\big(t_{\rm mix}^{(n)} + \alpha w_n\big) & = & 0  . 
\end{eqnarray}
For an introduction to Markov chains and mixing times see the 
%highly recommended 
book by Levin, Peres and Wilmer \cite{LevinPeresWilmer2009}, from which we borrow heavily.

\section{The Erd\H{o}s-R\'enyi random graph model}
The easiest model in the class of unimodal and radially symmetric distributions $\pi$ on the hypercube arises when we have $\theta_k = \theta \in (0,1]$ for all $k$, i.e. the acceptance probabilities for upward moves are constant across level sets. That is, $\pi(x) = \theta^{S(x)} (1+\theta)^{-n}$. The current paper will focus on this model.
The case where the neighboring odds ratios $\theta=\theta_k=\pi(v)/\pi(w)$ are bigger than one, where $S(v)=k+1$ and $S(w)=k$,  would correspond
to the mode of $\pi$ being at $\underline{1} = (1,1,...,1)$ instead of $\underline{0} = (0,...,0)$. By symmetry, this gives rise to nothing new, so we will assume $\theta \in (0,1]$ henceforth. The case where $\theta=1$ corresponds to $\pi$ being the uniform distribution. 

If we have $n = {\nu \choose 2}$ and identify the list of coordinates of the hypercube with the list of potential edges of a graph on $\nu$ vertices, then the hypercube represents the space of all possible (simple) graphs on $\nu$ vertices: a one indicates that a certain edge is present in the graph; a zero indicates that it is absent.
Since $\pi(x) = \left(\frac{\theta}{1+\theta}\right)^{\!S(x)} \left(1 - \frac{\theta}{1+\theta}\right)^{\!n-S(x)}$ for all $x \in \{0,1\}^n$, this distribution $\pi$ corresponds to the Erd\H{o}s-R\'enyi random graph model with parameter $p:=\frac{\theta}{1+\theta}$.
This is the probability distribution on (simple) graphs on $\nu$ vertices where each of the $n={\nu \choose 2}$ potential edges is present independently with probability $p$.

\vspace{0.5em}

The case where $\pi$ is the uniform distribution ($\theta=1$) corresponds to  the Erd\H{o}s-R\'enyi model with edge probability $\theta / (1+\theta) = 1/2$.
For this model it is well known that the (non-lazy) random walk Metropolis algorithm has cutoff at $(1/4) n \log n$ with a window of size $n$. 
This was proved by Aldous and Diaconis \cite{Aldous1983, DiaconisShahshahani1987} using Fourier analysis. 
The shape of the cutoff was studied in more detail by Diaconis, Graham and Morrison \cite{DiaconisGrahamMorrison1990}, whereas
Levin, Peres and Wilmer \cite[Theorem 18.3 on page 251]{LevinPeresWilmer2009} give a probabilistic proof.

The result has  been generalized to edge probabilities different from $1/2$: For general $\theta \in (0,1]$, it is known that
the non-lazy version of the random walk Metropolis chain 
has cutoff at $\frac{1}{2(1+\theta)} n \log n$ with window size $n$.
This was derived using Fourier analysis by Diaconis and Ram \cite[Theorem 5.4 on page 177]{DiaconisRam2000}.
See also the work of Diaconis and Hanlon \cite[Theorem 2 on page 104]{DiaconisHanlon1992} who explicitly calculate the eigenvalues and eigenvectors of the  transition kernel of the projection $S(X_t)$.
Ross and Xu \cite{RossXu1994} view the Metropolis chain for this model  as a random walk on a hypergroup deformation of the hypercube, 
and proceed by performing  Fourier analysis of this random walk.
See also \cite[page 2117]{DiaconisSaloff-Coste2006} for a discussion of these and some related results.
For the lazy version of this chain for  $\theta \in (0,1]$, we therefore get cutoff at $\frac{1}{1+\theta} n \log n$ with window size $n$.
This follows from 
%\cite[Theorem 5.4 on page 177]{DiaconisRam2000}
\cite{DiaconisRam2000}: For the upper bound, note that running the lazy chain for $2t$ steps corresponds to running the non-lazy chain for $T \sim Bin(2t,1/2)$ steps. So the standard deviation of $T$ is of order $\OO{\sqrt{t}} = \OO{\sqrt{n \log n}}$, which is easily absorbed into the window size term of order $n$. 
For the lower bound, note that making the chain lazy sends an eigenvalue $\lambda$ of the transition kernel to $(\lambda+1)/2$, while leaving the associated eigenfunction unchanged. Applying the results and methods of \cite[page 178]{DiaconisRam2000} then completes the proof.

In this paper we will give 
an alternative 
proof of this result, generalizing the methods used by Levin, Peres and Wilmer \cite[Theorem 18.3 on page 251]{LevinPeresWilmer2009} for the case where $\theta=1$. This is done in the hope that a probabilistic proof will be easier to generalize to less symmetric models, where 
Fourier analysis might be harder to apply. To be specific, we will proof the following result:

\begin{theorem}[{\bfseries Diaconis, Ram, 2000}]
\label{theta_model}
The lazy random walk Metropolis chain for $\pi(x) = \theta^{S(x)} (1+\theta)^{-n}$ on $\{0,1\}^n$ has cutoff at $\frac{1}{1+\theta} n \log n$ with a window of size $n$.
\end{theorem}

\begin{corollary}
Let $n := {\nu \choose 2}$ and let $\pi(x) = p^{S(x)} (1-p)^{n-S(x)}$ be the Erd\H{o}s-R\'enyi random graph model on
$\nu$ vertices with parameter $p \in (0,1)$. 
The lazy random walk Metropolis chain for this model  has cutoff at $\max\{p, 1-p\} n \log n$ with a window of size $n$.
\end{corollary}

\noindent {\bfseries Proof of the Corollary.} Let $\theta:=p/(1-p)$. As mentioned above, by relabeling states (switching zeros and ones) if necessary, we may assume $\theta \in (0,1]$. Since 
$$
p^{S(x)} (1-p)^{n-S(x)} = \left(\frac{\theta}{1+\theta}\right)^{\!S(x)} \left(1 - \frac{\theta}{1+\theta}\right)^{\!n-S(x)}
= \theta^{S(x)} (1+\theta)^{-n} 
$$
and $1/(1+\theta) = 1-p = \max\{p, 1-p\}$, the result follows from the Theorem. $\; \Box$

\section{Lower bound}
\label{sec_lower_bound_1}
Our proof for the lower bound part of Theorem \ref{theta_model} generalizes the proof given by Levin, Peres and Wilmer \cite[Proposition 7.13 on page 95]{LevinPeresWilmer2009} for the case where the stationary distribution $\pi$ is uniform ($\theta = 1$).
It is based on the method of distinguishing statistics, described as follows in \cite[Proposition 7.8 on page 92]{LevinPeresWilmer2009}: 

\begin{prop}[\bfseries{Levin, Peres, Wilmer, 2009}]
\label{dist_stat}
Let $\mu$ and $\nu$ be two probability distributions on $\mathcal{X}$, and let $S$ be a real-valued function on $\mathcal{X}$. If
$$
|E_{\mu}(S) - E_{\nu}(S)| \geq r \sigma ,
$$
where $\sigma^2 = [{\rm Var}_{\mu}(S) + {\rm Var}_{\nu}(S)]/2$, then
$$
||\mu - \nu||_{TV} \geq 1 - \frac{4}{4+r^2}  .
$$
\end{prop}

Here $E_{\mu}(S) := \sum_{x \in \mathcal{X}} S(x) \mu(x)$ denotes the expectation of $S$ under $\mu$, and likewise for $\nu$. So if we can find a real function on the state space $\mathcal{X}$
such that its expectations under $P^t(x,\cdot)$ and $\pi$ are still very different on the scale of the square root of their average variance after $t$ steps of the chain, then we have demonstrated that $||P^t(x,\cdot) - \pi||$ must still be large.

A natural choice for the distinguishing statistic is the number of ones in a state, $S(x) := \sum_{i=1}^n x_i$ for $x \in \{0,1\}^n$. 
Therefore we have to analyze the one-dimensional projection $S(X_t) =: S_t$ of our Markov chain $(X_t)$. As mentioned in the introduction, this is again a Markov chain whose transition probabilities satisfy 
$P(k,l) = P(x,S^{-1}(l)) = P(x,[y])$ 
for any $x,y$ with $S(x)=k$ and $S(y)=l$. As before, $[y] = \{z \in \{0,1\}^n \, : \, S(z)=S(y)\}$ denotes the equivalence class of all states with the same number of ones as $y$.
Because of their different domains it should not lead to confusion that we are using the same notation $P(\cdot,\cdot)$ for the transition probabilities of the original chain and the projected chain.

Similarly, the stationary distribution $\pi_S$ of $(S_t)$ is the push-forward $\pi_S := \pi S^{-1}$ of $\pi$ under $S$. This entails 
\begin{eqnarray*}
\pi_S(k)  & = & {n \choose k} \theta^k (1+\theta)^{-n} \\
& = & {n \choose k} \left( \frac{\theta}{1+\theta} \right)^{\!k} \left( 1 - \frac{\theta}{1+\theta} \right)^{\!n-k} \\
& = & \mbox{Binomial}\left(n, \frac{\theta}{1+\theta}\right)(k) .
\end{eqnarray*}
Therefore we get the expectation and variance of $S \sim \pi_S$ as
\begin{equation}
E_{\pi_S} S = \frac{n \; \theta}{1 + \theta}, \;\;\;\;\;
{\rm Var}_{\pi_S} S = \frac{n \; \theta}{(1 + \theta)^2}  .
\end{equation}

The chain $(S_t)$ is a birth and death chain on $\{0,1,...,n\}$ with transition probabilities
\begin{eqnarray}
\label{transition_kernel_1d}
 P(S_{t+1} = k+1 | S_t = k) & = & \left( 1 - \frac{k}{n} \right) \frac{\theta}{2}  , \nonumber \\
 P(S_{t+1} = k   | S_t = k) & = & \frac{1}{2} + \left( 1 - \frac{k}{n} \right) \frac{1-\theta}{2} ,  \\
 P(S_{t+1} = k-1 | S_t = k) & = & \frac{k}{2n}   . \nonumber
\end{eqnarray}
We begin by calculating its expectation after $t$ steps starting from $k$. Note that for all $t$ we get
$$
S_{t+1} - S_t = \left\{ \begin{array}{r@{\quad\mbox{with probability }\quad}l}
1 & \left( 1 - \frac{S_t}{n} \right) \frac{\theta}{2} , \\
-1 & \frac{S_t}{2n}  , \\
\end{array} \right.
$$ 
and $S_{t+1} - S_t = 0$ otherwise. So 
$$
E[S_{t+1} - S_t | S_t] = \left( 1 - \frac{S_t}{n} \right) \frac{\theta}{2} - \frac{S_t}{2n} =
\frac{\theta}{2} - S_t \frac{1+\theta}{2n}  ,
$$
and therefore
$$
E[S_{t+1} | S_t] = 
\frac{\theta}{2} + \left( 1 - \frac{1+\theta}{2n} \right) S_t  .
$$
By taking expectation $E_k$ with respect to the starting state $k$ we get for all $t,k$ that
\begin{equation}\label{exp3}
E_k (S_{t+1}) = \frac{\theta}{2} + \left( 1 - \frac{1+\theta}{2n} \right) E_k (S_t) .
\end{equation}
By induction on $t$ this leads to the following result:

\begin{prop}
\label{exp_theta}
The projected chain $S_t := S(X_t)$ of our lazy Metropolis chain
$(X_t)$ has for all $k=0,1,...,n$ and all $t \in \mathbb{N}$
\begin{equation}
\label{expect_St}
E_k (S_t) = \frac{n \theta}{1 + \theta} \left( 1 - \gamma^t \right) + k \gamma^t  ,
\end{equation}
where $\gamma := \gamma_{n,\theta} := 1 - \frac{1+\theta}{2n}$.
\end{prop}

%\emph{Proof: } By induction on $t$. Fix any starting state $k \in \{0,1,...,n\}$. The statement is true for $t=0$ since $E_k (V_0) = k$. Now suppose it is true for $t$. 
%and set $\gamma := \left( 1 - \frac{1+\theta}{2n} \right)$. 
%Then by (\ref{exp3}) and the induction hypothesis we get that
%\begin{eqnarray*}
%E_k (V_{t+1}) & = & \frac{\theta}{2} + \gamma \left[ \gamma^t \left( k - \frac{n \theta}{1+\theta} \right) + \frac{n \theta}{1+\theta} \right] \\
%& = & \frac{\theta}{2} + \gamma^{t+1} \left(  k - \frac{n \theta}{1+\theta} \right) + \left( \frac{n \theta}{1+\theta} - \frac{\theta}{2}  \right) \\
%& = & \gamma^{t+1} \left[ k - \frac{n \theta}{1+\theta} \right] + \frac{n \theta}{1+\theta}  ,
%\end{eqnarray*} 
%
%so it is also true for $t+1$. $\Box$

\begin{remark}
Note that the expected location  $E_k S_t$ is a convex combination of the starting state $S_0=k$ and the stationary mean  $E_{\pi_S} S$, with relative weights $\gamma^t$ and $1-\gamma^t$ respectively.
\end{remark}

\begin{remark}
\label{gamma^u}
By expanding $\log(\gamma)$ about one, it is easy to see that for $u := u_{n,\theta} := \frac{1}{1+\theta} n \log n$ we get
$$
\gamma^u \sim n^{-1/2}  ,
$$ 
by which we mean that $\lim_{n \rightarrow \infty} \frac{\gamma^u}{n^{-1/2}} = 1$.
\end{remark}

It remains to bound the variance of $S_t$. Since we want a lower bound on
$$
d(t) = \sup_x || P^t(x, \cdot) - \pi || \geq || P^t(\underline{1}, \cdot) - \pi ||  ,
$$ 
it's enough to consider the starting state $\underline{1} = (1,...,1)$ of all ones.
For this, first note that we can run the chain $X_t$ in the following way. For $t=0,1,2,...$, given we are at state $X_t$ at time $t$:
\begin{itemize}
\item Pick a coordinate $i \in [n]$ uniformly at random, independent of all previous choices.
\item Draw $U_t \sim$ Uniform$[0,1]$, independent of all previous choices.
\item Set $X_{t+1}^{(j)} := X_t^{(j)}$ for $j \ne i$, and set the $i^{th}$ coordinate of $X_{t+1}$ to
\end{itemize}
$$
X_{t+1}^{(i)} :=  \left\{ \begin{array}{l@{\quad:\quad}l}
1 & 0 \leq U_t \leq \frac{\theta}{2} , \\
X_t^{(i)} & \frac{\theta}{2} < U_t \leq \frac{1}{2} , \\
0 & \frac{1}{2} < U_t \leq 1   .\\
\end{array} \right.
$$
Now say that a coordinate $j$ has been \emph{refreshed} by time $t$, if coordinate $j$ was selected at some time $s<t$ \emph{and} $U_s \notin \left( \frac{\theta}{2}, \frac{1}{2} \right]$. 
Let $R_t$ be the number of coordinates not refreshed by time t. We can study the expectation and variance of $R_t$ with a natural modification of the classical coupon collector problem, where coordinate $j$ being refreshed corresponds to coupon $j$ being collected. This leads to the following result, with a proof analogous to the one for  \cite[Lemma 7.12 on page 94]{LevinPeresWilmer2009}:

\begin{prop}
\label{coupon_collector}
Consider the coupon collector problem with $n$ distinct coupon types, where at each trial, with probability $\frac{1-\theta}{2}$ we get \emph{no} coupon, and with probability $1 - \frac{1-\theta}{2} = \frac{1+\theta}{2}$ we get a coupon chosen (independently and) uniformly at random. Let $I_j(t)$ be the indicator of the event that the $j^{th}$ coupon has \emph{not} been collected by time $t$. Let $R_t := \sum_{j=1}^n I_j(t)$ be the number of coupon types not collected by time $t$. The random variables $I_j(t)$ are negatively correlated, and letting $\gamma :=  1 - \frac{1+\theta}{2n}$, we get for $t \geq 0$ that
\begin{eqnarray*}
E(R_t) & = & n  \gamma^t , \\
{\rm Var}(R_t) & \leq & n  \gamma^t \,(1-\gamma^t) \leq n  \gamma^t  . 
\end{eqnarray*}
\end{prop}

%
%
%
\begin{comment}
\vspace{0.5em}

{\emph Proof: } By definition of $I_j(t)$ we get 
$$
E[I_j(t)] = P\{\mbox{No coupon of type $j$ in trials 1,...,t}\} = \left( 1 - \frac{1+\theta}{2n} \right)^t = \gamma^t
$$ 
and
$$
{\rm Var}[I_j(t)] = E[(I_j(t))^2] - (E[I_j(t)])^2 = \gamma^t - \gamma^{2t} = \gamma^t \, (1-\gamma^t)  .
$$
Similarly, for $j \neq k$ we get
\begin{eqnarray*}
E[I_j(t) \, I_k(t)] & = & P\{\mbox{No coupon of type $j$ or $k$ in trials 1,...,t}\} \\
& = & \left( 1 - \frac{1+\theta}{2} \frac{2}{n} \right)^t \\
& = & \left( 1 - \frac{1+\theta}{n} \right)^t  , 
\end{eqnarray*}
so
$$
{\rm Cov}[I_j(t), I_k(t)] = E[I_j(t) \, I_k(t)] - E[I_j(t)] \, E[I_k(t)] = \left( 1 - \frac{1+\theta}{n} \right)^t - \left( 1 - \frac{1+\theta}{2n} \right)^{2t} \leq 0 .
$$
Therefore
$$
E[R_t] = \sum_{j=1}^n E[I_j(t)] = n \gamma^t
$$
and
$$
{\rm Var}[R_t] = \sum_{j=1}^n {\rm Var}[I_j(t)] + \sum_{j \neq k} {\rm Cov}[I_j(t), I_k(t)] \leq n \gamma^t \, (1-\gamma^t)  . \;\; \Box
$$

\end{comment}
%
%
%

\vspace{0.5em}

If we start the chain 
at $X_0 = \underline{1}$, then the conditional distribution of $S_t := S(X_t)$ given $R_t = r$ is the same as that of $r + B$, where $B \sim \mbox{Binomial}(n-r, \frac{\theta}{1+\theta})$. Therefore,
$$
E_{\underline{1}} [S_t \, | \, R_t] = R_t + (n - R_t) \frac{\theta}{1+\theta} = \frac{R_t + n \theta}{1+\theta}  ,
$$
so by taking expectation we get
$$
E_{\underline{1}} [S_t] = \frac{E [R_t] + n \theta}{1+\theta} = \frac{n \gamma^t + n \theta}{1+\theta}  
\left( = \frac{n \theta}{1+\theta} (1-\gamma^t) + n \gamma^t \right)  ,
$$
confirming our result (\ref{expect_St}) for general starting states $k$  for the special case $k=n$. Furthermore, since
$$
{\rm Var}_{\underline{1}} [S_t] = {\rm Var}\left[ E_{\underline{1}} (S_t \, | \, R_t) \right]
 + E\left[ {\rm Var}_{\underline{1}} (S_t \, | \, R_t) \right]  ,
$$
we get
\begin{eqnarray*}
{\rm Var}_{\underline{1}} [S_t] & = & {\rm Var}\left[ \frac{R_t + n \theta}{1+\theta} \right] + E\left[ {\rm Var} \Big(\mbox{Binomial}(n-R_t, \frac{\theta}{1+\theta})\Big) \right] \\
& = & \frac{1}{(1+\theta)^2} \, {\rm Var} [R_t] + \frac{1}{(1+\theta)^2} \, (n - E [R_t]) \, \theta \\
& \leq & \frac{1}{(1+\theta)^2} \left[ n \gamma^t + (n - n \gamma^t) \theta \right] \\
& \leq & \frac{n}{(1+\theta)^2} . 
\end{eqnarray*}
To apply Proposition \ref{dist_stat}, observe that
$$
\sigma^2 := \frac{{\rm Var}_{P^t(\underline{1}, \cdot)} S + {\rm Var}_{\pi} S}{2} \leq 
{\max\{{\rm Var}_{\underline{1}} S_t, {\rm Var}_{\pi_S} S\}} \leq \frac{n}{(1+\theta)^2}  .
$$
So for $\alpha>0, t := t_{n, \alpha} := \frac{1}{1+\theta} n \log n - \alpha n$ and $\gamma = 1 - \frac{1+\theta}{2n}$, we get for any fixed $\varepsilon>0$ and large $n$ that 
\begin{eqnarray*}
\left| E_{\underline{1}} S_t - E_{\pi} S \right| & = & \left| \frac{n (\theta + \gamma^t)}{1 + \theta} - \frac{n \theta}{1+\theta} \right| \\
& = & \frac{n}{1+\theta} \left( 1 - \frac{1+\theta}{2n} \right)^{\!t} \\
& \geq &  \sigma \sqrt{n} \left( 1 - \frac{1+\theta}{2n} \right)^{\!n \left( \frac{1}{1+\theta} \log n - \alpha \right)} \\
& \geq & \sigma \sqrt{n} (1-\varepsilon) \exp\left\{- \frac{1+\theta}{2} \left( \frac{1}{1+\theta} \log n - \alpha \right)\right\} \\
& = & \sigma (1-\varepsilon) \exp\left\{\alpha \frac{1+\theta}{2}\right\} \\
& =: & \sigma \, r_{\alpha}  .
\end{eqnarray*}
By Proposition \ref{dist_stat}, this means $d(t) \geq ||P^t(\underline{1},\cdot) - \pi|| \geq 1 - \frac{4}{4 + r_{\alpha}^2}$, and therefore
$$
\lim_{\alpha \rightarrow \infty} \liminf_{n \rightarrow \infty} \, d(t_{n, \alpha}) \geq
\lim_{\alpha \rightarrow \infty} 1 - \frac{4}{4 + r_\alpha^2} = 1  .
$$
This finishes the proof of the lower bound part of Theorem \ref{theta_model}. $\;\; \Box$

\vspace{0.5em}

\section{Lower bound, alternative proof}
\label{sec_azuma}
The previous proof for the lower bound part of Theorem \ref{theta_model} might be hard to generalize to distributions $\pi$ where $\theta_k$ is not constant in $k$. (Recall that $\theta_k := \pi(v)/\pi(w)$ for any $v,w$ such that $S(v)=k+1, S(w)=k$.) Therefore, we give here an alternative proof for this result. 
For this we use a modified version of the method of distinguishing statistics (Proposition \ref{dist_stat}), that avoids the need to  approximate the variance of the statistic $S$ under $P^t(x,\cdot)$. Instead, we use the fact that our chain only makes \emph{local moves} to argue that $S_t = S(X_t)$ must be concentrated about its mean under $P^t(x,\cdot)$. 
This will follow from Azuma's inequality, applied to the conditional expectation martingale formed by
$Y_i := E [S_t \, | \, S_{0:i}]$, for $i=0,1,...,t$, since we can show that  $(Y_i)$ has  bounded differences.
Here we write $S_{k:l} := (S_i)_{k \leq i \leq l}$.

%
%
%
\begin{comment}

\begin{theorem}[Azuma's Inequality]:
\label{azuma}
Let $(Y_i)_{i=0,1,...}$ be a martingale with bounded differences, i.e.
$$
a_i \leq Y_i - Y_{i-1} \leq b_i \qquad\mbox{for some } a_i, b_i \in \mathbb{R} \mbox{ and all } i \in \mathbb{N}  . 
$$
Then for any $t \in \mathbb{N}$ and $s>0$ we get
\begin{eqnarray*}
P\{Y_t - Y_0 \geq s\} & \leq & e^{-2s^2/c} \qquad\mbox{ and } \\
P\{Y_t - Y_0 \leq -s\} & \leq & e^{-2s^2/c}  ,
\end{eqnarray*}
where $c := \sum_{i=1}^t (b_i - a_i)^2$. 
\end{theorem}

For a proof see for example Dubhashi and Panconesi (2009, Theorem 5.2 on page 67), or Ross (1996, Theorem 6.3.3 on page 307). 

\end{comment}
%
%
%

\begin{prop}
\label{conc_theta}
Let $(X_t)$ be the lazy random walk Metropolis chain on $\{0,1\}^n$ for $\pi(x) = \theta^{S(x)} (1+\theta)^{-n}$. Let $S_t := S(X_t)$ and 
$v:=S_0=S(X_0)$. Then for any $s > 0$ and $\gamma := 1 - \frac{1+\theta}{2n}$ we get
\begin{eqnarray*}
P_v\{S_t \geq E_v S_t + s\} & \leq & e^{-\frac{2}{9}s^2(1-\gamma^2)}  \qquad\mbox{ and } \\
P_v\{S_t \leq E_v S_t - s\} & \leq & e^{-\frac{2}{9}s^2(1-\gamma^2)}  .
\end{eqnarray*}
\end{prop}

\vspace{0.5em}

\noindent \textbf{Proof. } Fix any $t \in \mathbb{N}$. Since $(S_i)$ is a Markov chain, we get from Proposition \ref{exp_theta} and the Markov property that for any $i=0,1,...,t$
$$
E[S_t \, | \, S_{0:i}] = E_{S_i} \, [S_{t-i}] = \gamma^{t-i} \left[ S_i - \frac{n \theta}{1+\theta} \right] + \frac{n \theta}{1+\theta}  .
$$
Therefore for any $i=1,2,...,t$, we get
\begin{eqnarray*}
E[S_t \, | \, S_{0:i}] - E[S_t \, | \, S_{0:(i-1)}] & = & \gamma^{t-i} \left[ S_i - \frac{n \theta}{1+\theta} \right] - \gamma^{t-(i-1)} \left[ S_{i-1} - \frac{n \theta}{1+\theta} \right] \\
 & = & \gamma^{t-i} \left[ S_i - S_{i-1} + (1-\gamma) S_{i-1} - (1-\gamma) \frac{n \theta}{1+\theta} \right] \\
 & = & \gamma^{t-i} \left[ S_i - S_{i-1} + \frac{1+\theta}{2n} S_{i-1} - \frac{\theta}{2} \right]  .
\end{eqnarray*}
This means the martingale $(Y_i)$ has bounded differences  
$$
| E[S_t \, | \, S_{0:i}] - E[S_t \, | \, S_{0:(i-1)}] | \leq \frac{3}{2} \gamma^{t-i}  ,
$$
since we have $S_i \in [0,n]$ and $|S_i - S_{i-1}| \leq 1$, and also $\theta \in (0,1]$.
Applying Azuma's inequality
now gives the result, since  we have
$$
c := 4 \sum_{i=1}^t \left(\frac{3}{2} \gamma ^{t-i}\right)^{\!2} = 9 \sum_{i=0}^{t-1} \gamma^{2i} \leq \frac{9}{1-\gamma^2} . \Box
$$

\vspace{0.5em}

Now we use this concentration result for $(S_t)$ 
to establish the lower bound on the mixing time for this Markov chain. 
Fix $\alpha > 0$ and set $t := t_{n,\alpha} :=  \frac{1}{1+\theta} n \log n - \alpha n$. 
We need to show that
\begin{equation}
\label{theta_lbb}
\lim_{\alpha \rightarrow \infty} \liminf_{n \rightarrow \infty} d(t_{n,\alpha}) = 1  .
\end{equation}
%(Note that since $d(t)$ is decreasing in $t$, it is enough to prove (\ref{theta_lb}) for \emph{large} $C<1$, so our restriction $\delta < \frac{1}{2}$ is not problematic.) 
Denote with $S_t := S(X_t)$ the projection of our chain under $S$, started at $X_0 = \underline{\mathbf{1}}$, so that $S_0=n$  and note that for any $r \in [0,n]$ we get 
\begin{eqnarray*}
d(t) & \geq & ||P^t(\underline{1},\cdot) - \pi|| \\
% 		 & =    & sup_{A \subset \{0,1\}^n} P^t(\underline{1},A) - \pi(A) \\
 		 & \geq & sup_{L \subset [n]} P^t(\underline{1},S^{-1}(L)) - \pi(S^{-1}(L)) \\
% 		 & =    & sup_{L \subset [n]} P^t(n,L) - \pi_S(L) \\
%     & \geq & P_n\{V_t \geq r\} - P\{B_n \geq r\} \\
     & \geq & 1 - P_n\{S_t < r\} - P\{B_n \geq r\}  .
\end{eqnarray*}
Here the second inequality comes from noting that the  supremum on the right is over a subset of events from the  supremum defining the total variation distance on the left. Therefore, total variation distance can only decrease after projecting down. 
For the last inequality we took the event $L := \{\lceil r \rceil ,\lceil r \rceil +1,...,n\}$ and let
$B_n$ stand for a random variable with distribution $\pi_S=\mbox{Binomial}(n,\frac{\theta}{1+\theta})$. 
Our task is to find  $r \in [0,n]$ such that 
\begin{eqnarray}
\label{conc_th}
\lim_{\alpha \rightarrow \infty} \limsup_{n \rightarrow \infty} P_n\{S_t < r\} & = & 0 \qquad\mbox{ and } \\
\label{conc_Bin}
\lim_{\alpha \rightarrow \infty} \limsup_{n \rightarrow \infty} P\{B_n \geq r\} & = & 0 . 
\end{eqnarray}
Suppressing rounding issues from the notation, we can pick 
$r := np + \sqrt{\alpha n}$, where we write $p := \frac{\theta}{1+\theta}$.
To establish (\ref{conc_th}), we use the concentration result for $S_t$ from Proposition \ref{conc_theta}. 
Using the results from Proposition \ref{exp_theta} and Remark \ref{gamma^u}, this shows that for all $\epsilon'>0$ and large $n$ we get
\begin{eqnarray*}
P_n\{S_t < r\}& = & P_n\{S_t - E_n S_t < -((1-p)n \gamma^t - \sqrt{\alpha n})\} \\
& \leq & P_n\left\{S_t - E_n S_t < -\left((1-\epsilon')(1-p)\sqrt{n} e^{\alpha (1+\theta)/2} - \sqrt{\alpha n}\right)\right\} \\
& \leq & \exp\left\{- \frac{2}{9} n \left( (1-\epsilon')(1-p) e^{\alpha (1+\theta)/2} - \sqrt{\alpha} \right)^2 (1-\gamma^2)\right\} \\
& \leq & \exp\left\{- \frac{2}{9} \left( (1-\epsilon')(1-p) e^{\alpha (1+\theta)/2} - \sqrt{\alpha} \right)^2 (1+\theta) \gamma\right\} .
\end{eqnarray*}
Here we use the fact that $1-\gamma^2 \geq \frac{1+\theta}{n} \gamma$ for the last inequality.
Taking limits
establishes (\ref{conc_th}).
The result (\ref{conc_Bin}) follows from Chebychev's inequality applied to the Binomial($n,p$) random variable $B_n$. $\Box$

\section{Upper bound}
To establish the upper bound part of Theorem \ref{theta_model}, 
we use a two-stage coupling procedure, generalizing  the one given by Levin, Peres and Wilmer \cite[Theorem 18.3 on page 251]{LevinPeresWilmer2009} for the case where $\pi$ is uniform ($\theta=1$). 
In that case, it is enough to study the one-dimensional projection $X_t \mapsto S(X_t)$, since due to symmetry, total variation distance to stationarity is the same from any starting state. Therefore, we may start our chain at the state of all zeros (or all ones), in which case total variation distance to stationarity doesn't change under the projection to the number of ones.

For general $\theta \in (0,1]$, total variation distance to stationarity is not necessarily the same from any starting state, and we weren't able to prove that the state of all zeros (or all ones) is a worst starting state. So it's not clear to us whether the problem of upper bounding the mixing time can still be reduced to a one-dimensional problem here. 
However, we were able to get a sharp upper bound on the mixing time using a \emph{two}-dimensional projection that depends on the starting state $X_0=x$. For this, consider $Z := Z_x \colon \{0,1\}^n \rightarrow \{0,1,...,n\} \times \{0,1,...,n\}$,
$$
X_t \mapsto Z_x(X_t) := Z_t := \left( S(X_t), d(x,X_t) \right)  ,
$$
where 
$d(x,y)$ is the number of coordinates where $x$ and $y$ disagree,
 and $x=X_0$ is the starting state of the chain. In words, we project down to the number of ones, $S(X_t)$, and the distance to the starting state, $d(X_t,X_0)$. This two-dimensional process $(Z_t)$ is similar to the two-coordinate chain used in
\cite[section 3]{LevinLuczakPeres2010} in the study of  Glauber dynamics for the mean-field Ising model.

Our proof proceeds by showing that $(Z_t)$ is a (two-dimensional) birth and death chain with the same total variation distance to its stationary distribution as the original chain $(X_t)$. Bounding this distance is then  achieved by a two-stage coupling procedure using an independence coupling of two versions of the chain. The first stage (\emph{drift}-regime) brings them close together in expectation after $\frac{1}{1+\theta} n \log n$ steps due to the drift towards the mean in this birth and death chain. In the second stage (\emph{entropy}-regime)  
the drift of $(Z_t)$ is likely weak, so a comparison to the behavior of simple random walk shows that
the chains can be made to coalesce after an additional $\alpha n$ steps with high probability.

\subsection{Properties of the projection}
We begin by establishing some properties of the two-dimensional projection $Z_t := Z_x(X_t)$ of $X_t$, for which we need some more notation. Fix $x,z \in \{0,1\}^n$ and let $S(x)=k, S(z)=l, d(x,z)=l'$. Define
\begin{eqnarray*}
G & := & G(x,z) := \{i \in [n] \, : \, x^{(i)}=0, z^{(i)}=0\} , \\
N & := & N(x,z) := \{i \in [n] \, : \, x^{(i)}=0, z^{(i)}=1\} , \\
E & := & E(x,z) := \{i \in [n] \, : \, x^{(i)}=1, z^{(i)}=0\} , \\
F & := & F(x,z) := \{i \in [n] \, : \, x^{(i)}=1, z^{(i)}=1\}  .
\end{eqnarray*}
When clear from the context, we might suppress the dependence of $G, N, E, F$ on $x, z$ in the notation. For the number of elements in these four sets we get
\begin{eqnarray}
\label{nG}
\#G & = &n - \frac{l + l' + k}{2} , \nonumber \\
\#N & = & \frac{l + l' - k}{2} , \\
\#E & = &\frac{l' + k - l}{2} , \nonumber  \\
\#F & = & \frac{l - (l' - k)}{2} \nonumber   .
\end{eqnarray}
%Here the last equality follows from
%$$
%l' - l + F = g - N = E = k - F \; .
%$$
%Together with the first equality in each of rows 1,2,3 above, this implies the second equality in rows 1,2,3. 
This follows from $\#N+\#E=l', \#N+\#F=l, \#E+\#F=k, \#G+\#N=n-k$, e.g. by starting with the observation 
$$
l' - l + \#F = l' - \#N = \#E = k - \#F  ,
$$
which gives the last equality $2 \#F  =  l - (l' - k)$ above. The remaining equalities then follow.
This is probably best understood by looking at an example. Suppose
\begin{eqnarray*}
x & = & 0 0 0 0 0 \, 0 0 0 \, 1 1 \, 1  ,\\
z & = & 0 0 0 0 0 \, 1 1 1 \, 0 0 \, 1 .
\end{eqnarray*}
Then $n=11, S(x)=k=3, S(z)=l=4, d(x,z)=l'=5$ and $\#G=5, \#N=3, \#E=2, \#F=1$.
In general, we have a one-to-one correspondence between $(n,k,l,l')$ and 
$(\#G,\#N,\#E,\#F)$. Note that the formulas for $\#G,\#N,\#E,\#F$ above all give integers because $l, l'-k, l'+k$ always have the same parity.

\vspace{0.5em}

\begin{prop}
\label{const_2d}
For each $x \in \{0,1\}^n$ the $t$-step transition probabilities $P^t(x,\cdot)$ are constant on the level sets 
$$
\mathcal{X}(l,l',x) := \left\{ z \in \{0,1\}^n \, : \, Z_x(z)=(l,l') \right\}
$$
for all $l,l' \in \{0,1,...,n\}$.
\end{prop}

\vspace{0.5em}

\noindent {\bfseries Proof.}  Let a permutation $\phi$ act on $\{0,1\}^n$ by mapping
$x=(x_1,...,x_n)$ to $\phi(x)=(x_{\phi(1)},...,x_{\phi(n)})$.
Fix $x \in \{0,1\}^n$. By symmetry, we get
\begin{eqnarray*}
P(x,y) & = & P\big( \phi(x),\phi(y)\big) , \\
P^t(x,y) & = & P^t\big(\phi(x),\phi(y)\big) , \\
P^t(x,A) & = & P^t\big(\phi(x),\phi(A)\big) , 
\end{eqnarray*}
for all states $y$, all times $t$, and all subsets $A$ of $\{0,1\}^n$. Now fix $l,l' \in \{0,1,...,n\}$ and pick any $y,z \in \mathcal{X}(l,l',x)$. Then there exists a permutation $\phi$ that maps $y$ to $z$ and leaves $x$ fixed; 
$\phi$ maps $\{i  :  x_i = 0\}$ to itself and $\{i  :  x_i = 1\}$ to itself. Thus we get
$$
P^t(x,y) = P^t\big(\phi(x),\phi(y)\big) = P^t(x,z)
$$
as desired. $\Box$

\vspace{0.5em}

\begin{corollary}
\label{proj_2d}
The projection $Z_t := Z_x(X_t) := \big(S(X_t), d(x,X_t)\big)$ is Markov and we get
$$
|| P^t(x,\cdot) - \pi || = || \mathcal{D}(_kZ_t) - \pi_{Z_x} ||  ,
$$
where $x=X_0$ and $k=S(x)$. Here, $\mathcal{D}(_kZ_t)$ denotes the distribution of the chain $(Z_t)$ at time $t$ when started at $Z_0 = (k,0)$ and 
$\pi_{Z_x} := \pi Z_x^{-1}$ is the stationary distribution of the chain $(Z_t)$.
\end{corollary}

\vspace{0.5em}

\noindent {\bfseries Proof.}  Fix $X_0 = x$ with $S(x)=k$ and consider the equivalence relation $\sim_x$ corresponding to the classes 
$\mathcal{X}(l,l',x)$,
for $l,l' \in \{0,1,...,n\}$. That is, for $y,z \in \{0,1\}^n$ we have
$$
y \sim_x z \; \Leftrightarrow \; Z_x(y) = Z_x(z)  . 
$$
Then the projection $(Z_t)$ of $(X_t)$ is a Markov chain, if
\begin{equation}
\label{proj_markov}
P\big(y,\mathcal{X}(l,l',x)\big) = P\big(z,\mathcal{X}(l,l',x)\big)
\end{equation}
for all $l,l' \in \{0,1,...,n\}$ and all $y,z \in \{0,1\}^n$ such that $y \sim_x z$. 
So fix any $y \sim_x z$ and any $l,l' \in \{0,1,...,n\}$. Then there exists a permutation $\phi$ of $1 \colon \!n$ such that $\phi(y)=z$ and $\phi(\mathcal{X}(l,l',x)) = \mathcal{X}(l,l',x)$. 
%
% Proof:
% Let \phi be the permutation mapping y to z from before (\phi only permutes within the (i:x_i=0) and (i:x_i=1).)
% Fix v \sim_x w \in \mathcal{X}(l,l',x)
% To show [v] \subset \phi([v]), let \psi be the permutation s.t. \psi(v)=w, \psi only permuting within (i:x_i=0) etc as above.
% Then let s := \phi^{-1}(\psi(v)). 
% we get s \sim_x v and \phi(s)=\psi(v)=w \Box
%
% Since \phi is a bijection, [v] \subset \phi([v]) implies [v] = \phi([v]). \Box \Box
%
Consequently, we get
$$
P\big(y,\mathcal{X}(l,l',x)\big) = P\big(\phi(y),\phi(\mathcal{X}(l,l',x))\big) = P\big(z,\mathcal{X}(l,l',x)\big) ,
$$
showing that $(Z_t)$ is Markov.

Total variation distance to stationarity remains unchanged under the projection $Z_x$ because both $P^t(x,\cdot)$ and $\pi$ are constant on sets 
$\mathcal{X}(l,l',x)$
for $l,l' \in \{0,1,...,n\}$ by Proposition \ref{const_2d}. That allows us to pull out the absolute values from the inner sum in the second equation below, because all the terms in the sum are equal:
\begin{eqnarray*}
|| P^t(x,\cdot) - \pi || & = & \frac{1}{2} \sum_l \sum_{l'} \sum_{z \in \mathcal{X}(l,l',x)} \left| P^t(x,z) - \pi(z) \right| \\
& = & \frac{1}{2} \sum_l \sum_{l'} \left| \sum_{z \in \mathcal{X}(l,l',x)}  P^t(x,z) - \pi(z) \right| \\
& = & \frac{1}{2} \sum_l \sum_{l'} \left| P^t\big(x, Z_x^{-1}(l,l')\big) - \pi\big(Z_x^{-1}(l,l')\big)  \right| \\
& = & || P^t(x,\cdot) Z_x^{-1} - \pi Z_x^{-1} ||  .
\end{eqnarray*}
Clearly, $P^t(x, \cdot) Z_x^{-1} = \mathcal{D}(_kZ_t)$, and the fact that $\pi_{Z_x} := \pi Z_x^{-1}$ is stationary for $(Z_t)$ is an elementary calculation. $\; \Box$

\subsection{Reparametrization}
For any fixed $X_0=x$ with $S(x)=k$ the projected chain $(Z_t) = (Z_x(X_t))$ has the following transition kernel:
\begin{equation} 
\label{transition_kernel_2d}
P((l,l'),(h,h')) =  \left\{ \begin{array}{l@{\quad:\quad}l}
\frac{2n - (l'+l+k)}{4n} \theta & h=l+1, h'=l'+1 , \\
\frac{l'+l-k}{4n} & h=l-1, h'=l'-1 , \\
\frac{k+l'-l}{4n} \theta & h=l+1, h'=l'-1 , \\
\frac{k-(l'-l)}{4n} & h=l-1, h'=l'+1 , \\
\frac{1}{2} + \frac{n-l}{2n} (1-\theta) & h=l, \;\;\;\;\;\; h'=l' , \\
0 & \mbox{otherwise}  . \\  
\end{array} \right.
\end{equation}
This follows from (\ref{nG}) together with the transition rule of the original chain $(X_t)$. For example,
$$
P((l,l'),(l+1,l'+1)) = P\{ \mbox{pick } i \in G \mbox{ and flip } 0 \rightarrow 1\} = \#G \frac{1}{2n} \theta  .
$$
For $k \leq \frac{n}{2}$ the state space of this chain is
$$
\begin{array}{ll@{\quad:\quad}l}
\{(l,l') \in \{0,1,...,n\}^2 \, : \, & 
l' \in k + \{-l, -l+2, ..., l\} & \mbox{ for } l \leq k , \\
& l' \in k + \{l-2k, l-2k+2,...,l\} & \mbox{ for } k < l < n-k  , \\
& l' \in n-k + \{l-n, l-n+2,...,n-l\} & \mbox{ for } n-k \leq l  
\}  .\\
\end{array}
$$
A similar result holds for $k \geq \frac{n}{2}$.
% That is,
%$$
%\begin{array}{ll@{\quad:\quad}l}
%\{(l,l') \in \{0,1,...,n\}^2 \, : \, & 
%l' \in \{k-l, k-l+2, ..., k+l\} & \mbox{ for } l \leq k\; , \\
%& l' \in \{l-k, l-k+2, ..., l+k\} & \mbox{ for } k < l < n-k \; , \\
%& l' \in \{l-k, l-k+2, ..., 2n-(k+l)\} & \mbox{ for } n-k \leq l  
%\} \; .\\
%\end{array}
%$$
%
%Alternative formulation:
%$$
%\begin{array}{ll@{\quad:\quad}l}
%\{(l,l') \in \{0,1,...,n\}^2 \, : \, & 
%l' \in k + \{-l, -l+2, ..., l\} & \mbox{ for } l \leq k\; , \\
%& l' \in k + \{-l+2(l-k), -l+2(l-k)+2,...,l\} & \mbox{ for } k < l < n-k \; , \\
%%% or, equivalently:
%%%& l' \in n-k + \{l-n, l-n+2,...,n-l-2(n-k-l)\} & \mbox{ for } k < l < n-k \; ,  \\
%& l' \in n-k + \{l-n, l-n+2,...,n-l\} & \mbox{ for } n-k \leq l  
%\} \; .\\
%\end{array}
%$$
In both cases, after reparametrizing
$$
(l,l') \mapsto (l'-l, l'+l) =: (r,r') ,
$$
the state space becomes $\{(r,r') \, : \, r \in \{-k,-k+2,...,k\}, r' \in \{k,k+2,...,2n-k\}$.
%, isomorphic to a rectangle in $\mathbb{Z}^2$
The boundaries $-k \leq r \leq k$ and $k \leq r' \leq 2n-k$ here can also be confirmed like this:
By definition we have 
\begin{eqnarray*}
r  & = & l'-l = N+E-(N+F) = E-F \;\;\;\;\;\; \mbox{ and} \\ 
r' & = & l'+l = N+E+(N+F) = 2N+k . 
\end{eqnarray*}
Since $F \geq 0$, we get $r=E-F \leq E \leq E+F = k$. Also, since $E \geq 0$ and $F \leq k$, we get $r=E-F \geq -F \geq -k$. Similarly, since $N \leq n-k$, we get $r'=2N+k \leq 2(n-k)+k=2n-k$. And since $N \geq 0$, we get $r'=2N+k \geq k$.

The transition kernel in this new parametrization becomes
\begin{equation} 
\label{transition_kernel_2dr}
P((r,r'),(s,s')) =  \left\{ \begin{array}{l@{\quad:\quad}l}
\frac{2n - (r'+k)}{4n} \theta & s=r,  \;\;\;\;\;\;\, s'=r'+2 , \\
\frac{r'-k}{4n} & s=r,  \;\;\;\;\;\;\, s'=r'-2 , \\
\frac{k+r}{4n} \theta & s=r-2, \, s'=r' , \\
\frac{k-r}{4n} & s=r+2, \, s'=r' , \\
\frac{1}{2} + \frac{2n-(r'-r)}{4n} (1-\theta) & s=r, \;\;\;\;\;\;\,\, s'=r' , \\
0 & \mbox{otherwise}  .\\  
\end{array} \right.
\end{equation}
So the chain $(Z_t)$ can be viewed as a birth and death chain on a rectangle in $\mathbb{Z}^2$. A useful feature of the parametrization (\ref{transition_kernel_2dr}) is that here the probability of moving up (down) in the $r$-dimension only depends on the current location in that dimension: it only depends on $r$, not on $r'$. Similarly, the probability of moving up (down) in the $r'$-dimension only depends on $r'$, not on $r$. Therefore, the problem of coupling two versions of this chain can  be split up into coupling two one-dimensional processes. For this reason we will use the parametrization (\ref{transition_kernel_2dr}) for the rest of the paper.
%Note also that this feature will be lost if we allow $\theta=\theta_{S(X_t)}$ to depend on the number of ones in the current state $X_t$.
Note that the chain $(Z_t)$ and its stationary distribution depend on the initial state $x$ (and its number of ones $k$) used for the projection.

\subsection{Expected location of $(Z_t)$}
We now calculate the expected location of the chain $(Z_t)$ after $t$ steps when started at $Z_0 = (r,r')$ . Similar to the one-dimensional projection $S(X_t)$ that we analyzed for the lower bound, this expectation can be calculated explicitly by induction on $t$, since the transition probabilities (\ref{transition_kernel_2dr}) are all linear in the current location $(r,r')$ of the chain.

Fix $X_0=x \in \{0,1\}^n$ and let $k=S(x), Z_t=Z_x(X_t)=(S(X_t),d(x,X_t))$, so that $Z_0=(k,0)$. 
Denote with $Z_t = (Z_t^{(r)},Z_t^{(r')})$ the coordinates of the chain in the new parametrization (\ref{transition_kernel_2dr}), so that $(Z_0^{(r)},Z_0^{(r')}) = (-k,k)$, and write $E_k$ for the expectation operator given this starting state. 
Then for this parametrization,
$$
Z_{t+1} - Z_t =  \left\{ \begin{array}{l@{\quad:\quad}l}
(0,2)   & \mbox{ with probability } \;\; \frac{2n - (k+Z_t^{(r')})}{4n} \, \theta  , \\  
(0,-2)  & \mbox{ with probability } \;\; \frac{Z_t^{(r')}-k}{4n}   , \\  
(-2,0)  & \mbox{ with probability } \;\; \frac{k+Z_t^{(r)}}{4n} \, \theta  , \\  
(2,0)   & \mbox{ with probability } \;\; \frac{k-Z_t^{(r)}}{4n}   , \\  
(0,0)   & \mbox{ otherwise }  . \\  
\end{array} \right.
$$
Therefore,
\begin{eqnarray}
\label{drift_function}
& & \hspace{-3em} E[Z_{t+1}-Z_t \, | \, Z_t] \nonumber\\
& = & \left( -2 \, \frac{k+Z_t^{(r)}}{4n} \, \theta + 2 \, \frac{k-Z_t^{(r)}}{4n} \, , \,
2 \, \frac{2n - (k+Z_t^{(r')})}{4n} \, \theta - 2 \frac{Z_t^{(r')}-k}{4n} \right)  \nonumber \\
& = & \left( \frac{k(1-\theta)-Z_t^{(r)}(1+\theta)}{2n}, \frac{2n\theta + k(1-\theta)-Z_t^{(r')}(1+\theta)}{2n} \right) , \end{eqnarray}
so that $E[Z_{t+1} \, | \, Z_t]$ is equal to 
$$  
\left( \frac{k(1-\theta)+[2n - (1+\theta)] Z_t^{(r)}}{2n}, \frac{2n\theta + k(1-\theta) + [2n - (1+\theta)] Z_t^{(r')}}{2n} \right)  .
$$
By taking expectation, we get 
\begin{eqnarray}
\label{indstep2}
& &  \hspace{-2em} E_k [Z_{t+1}] \nonumber \\
& = & \left( \frac{k(1-\theta)+[2n - (1+\theta)] E_k Z_t^{(r)}}{2n}, \frac{2n\theta + k(1-\theta) + [2n - (1+\theta)] E_k Z_t^{(r')}}{2n} \right) \; \nonumber \\
&  =: &  \left( \beta + \gamma \, E_k Z_t^{(r)}, \theta + \beta + \gamma \, E_k Z_t^{(r')} \right)  .
\end{eqnarray}
%where $\beta := \beta_{n,k,\theta} := \frac{k}{2n} (1-\theta)$ and $\gamma := \gamma_{n,\theta} := 1 - \frac{1+\theta}{2n}$.
By induction on $t$, this leads to a proof of the following result:

\begin{prop}
\label{exp_2d}
Let $Z_t = Z_x(X_t) = \big(S(X_t),d(x,X_t)\big)$ be the two-dimensional projection of the lazy random walk Metropolis chain $(X_t)$ for $\pi(x) = \theta^{S(x)} (1+\theta)^{-n}$, started at $X_0 = x \in \{0,1\}^n$ with $S(x)=k$.
Then, in the parametrization (\ref{transition_kernel_2dr}) and for any $t \in \mathbb{N}$, we get
$$
E_k \left(Z_t^{(r)},Z_t^{(r')}\right) = \left( \frac{2n\beta}{1+\theta} (1-\gamma^t) - k \gamma^t \, , \,
\frac{2n (\theta + \beta)}{1+\theta} (1 - \gamma^t) + k \gamma^t \right),
$$
where $\beta := \beta_{n,k,\theta} := \frac{k}{2n} (1-\theta)$ and $\gamma := \gamma_{n,\theta} := 1 - \frac{1+\theta}{2n}$.
\end{prop}

\vspace{0.5em}

\noindent {\bfseries Proof.} The claim is true for $t=0$.
Now suppose it is true for $t$. Then by (\ref{indstep2}) we get that $E_k \left(Z_{t+1}^{(r)}, Z_{t+1}^{(r')}\right)$ is equal to
\begin{eqnarray*}
 &  & \left( \beta + \gamma \left[ \frac{2n\beta}{1+\theta} (1-\gamma^t) - k \gamma^t \right] \, , \,  \theta + \beta + \gamma \left[ \frac{2n(\theta+\beta)}{1+\theta} (1-\gamma^t) + k \gamma^t \right] \right) \\
& = & \left( \beta + \frac{2n\beta}{1+\theta} \gamma - \frac{2n\beta}{1+\theta} \gamma^{t+1} - k \gamma^{t+1} \, , \, 
\theta + \beta + \frac{2n(\theta+\beta)}{1+\theta} \gamma - \frac{2n(\theta+\beta)}{1+\theta} \gamma^{t+1} + k \gamma^{t+1} \right) \\
& = & \left( \frac{2n\beta}{1+\theta} \left[ \frac{1+\theta}{2n} + \gamma \right] - \frac{2n\beta}{1+\theta} \gamma^{t+1} - k \gamma^{t+1} \, , \, 
\frac{2n(\theta+\beta)}{1+\theta} \left[ \frac{1+\theta}{2n} + \gamma \right] - \frac{2n(\theta+\beta)}{1+\theta} \gamma^{t+1} + k \gamma^{t+1} \right) \\
& = & \left( \frac{2n\beta}{1+\theta} (1-\gamma^{t+1}) - k \gamma^{t+1} \, , \, 
\frac{2n(\theta+\beta)}{1+\theta} (1-\gamma^{t+1}) + k \gamma^{t+1} \right)  .
\end{eqnarray*}
So the claim is also true for $t+1$. $\; \Box$

\vspace{0.5em}

\begin{corollary}
\label{stat_exp_2dr}
For the expectation under stationarity, in the parametrization (\ref{transition_kernel_2dr}) of our two-dimensional chain $(Z_t)$, we get
$$
E_{\pi} \left( Z^{(r)}, Z^{(r')} \right) = \left( \frac{2n\beta}{1+\theta}\, , \, \frac{2n(\theta+\beta)}{1+\theta} \right)
= \left( k \frac{1-\theta}{1+\theta} \, , \, \frac{2n\theta}{1+\theta} + k \frac{1-\theta}{1+\theta} \right)  .
$$
\end{corollary}

\vspace{0.5em}

\noindent {\bfseries Proof.} Since the Markov chain $(Z_t)$ is irreducible and aperiodic, it converges to its unique stationary distribution as $t$ goes to infinity for fixed $n,k$. Since the state space is finite, this convergence also holds 
for expectations.
So by the Proposition, we get
\begin{eqnarray*}
E_{\pi} \left( Z^{(r)}, Z^{(r')} \right) & = & \lim_{t \rightarrow \infty} E_k \left( Z_t^{(r)},Z_t^{(r')} \right) \\
& = & \left( \frac{2n\beta}{1+\theta} \, , \, \frac{2n(\theta+\beta)}{1+\theta} \right)  ,
\end{eqnarray*}
since for fixed $n,k$ we have $\gamma^t \rightarrow 0$ as $t$ goes to infinity. $\;\;\; \Box$

\vspace{0.5em}

\begin{remark}
By reversing the linear transformation $(l,l') \mapsto (l'-l,l'+l)$, we immediately get
$$
E_{\pi} Z = \left( \frac{n\theta}{1+\theta} \, , \, \frac{n\theta + k(1-\theta)}{1+\theta} \right)
$$
for the expectation under stationarity in the original parametrization (\ref{transition_kernel_2d}) of our two-dimensional chain $(Z_t)$.
For the first coordinate this confirms what we already know from $S(X) \sim \mbox{Binomial}(n, \frac{\theta}{1+\theta})$ under stationarity. 
\end{remark}

\begin{remark}
\label{exp_2dr_general}
By exactly the same proof we get for a general starting state $(v,v')$ in the new parametrization (\ref{transition_kernel_2dr}) 
$$
E_{(v,v')} \left(Z_t^{(r)},Z_t^{(r')}\right) = \left( \frac{2n\beta}{1+\theta} (1-\gamma^t) + v \, \gamma^t \, , \,
\frac{2n (\theta + \beta)}{1+\theta} (1 - \gamma^t) + v' \gamma^t \right)  .
$$
Note that the expected location of $Z_t$ after $t$ steps is a convex combination of the starting state $Z_0=(v,v')$ and the stationary mean $E_{\pi_{Z_x}} Z$, with relative weights $\gamma^t$ and $1-\gamma^t$ respectively.
\end{remark}

\begin{remark}
\label{end_stage_one}
From Remark \ref{gamma^u} on page \pageref{gamma^u}, we know that for $u := u_{n,\theta} := \frac{1}{1+\theta} n \log n$ we get $\gamma^u \sim n^{-1/2}$.
Therefore, Proposition \ref{exp_2d} and Corollary \ref{stat_exp_2dr} imply that for any starting state the expected location of the chain $(Z_t)$ after $u$ steps is within $O(\sqrt{n}\,)$ of the expected location of the chain under stationarity. 
\end{remark}

\vspace{0.5em}

We now want to show that an additional $\alpha n$ number of steps is enough to couple (with high probability) two chains that are at distance $O(\sqrt{n}\,)$ of their stationary mean. This will follow from a corresponding result for simple random walk on $\mathbb{Z}^2$, since close to the stationary mean we are now in the ``entropy regime'' where the drift of the chain $(Z_t)$ is negligible.

\subsection{Burn-in}
But first we show that by running the chain for an initial $\alpha n$ steps (burn-in period), we may assume that the number of ones in the state that is used in the two-dimensional projection is close to its stationary mean. Fix $\delta>0$ and 
let $p := \frac{\theta}{1+\theta}$. 
We will show that
\begin{equation}
\label{burnin_1}
\max_{x \in \{0,1\}^n} || P^{\alpha n + t}(x,\cdot) - \pi || \leq \max_{y: S(y) \in n (p \pm \delta)} || P^{t}(y,\cdot) - \pi || + o(1) ,
\end{equation}
where the $o(1)$ term goes to zero as $n$ goes to infinity (uniformly in $t$).
To see this, we condition on where we are after the first $\alpha n$ steps:
\begin{eqnarray*}
& &  || P^{\alpha n + t}(x,\cdot) - \pi || \\
& = & || \sum_{y} P^{\alpha n}(x,y) \left[ P^{t}(y,\cdot) - \pi \right] || \\
& \leq &  \sum_{y: S(y) \in    n(p \pm \delta)} P^{\alpha n}(x,y) \; || P^{t}(y,\cdot) - \pi || +
          \sum_{y: S(y) \notin n(p \pm \delta)} P^{\alpha n}(x,y) \; || P^{t}(y,\cdot) - \pi || \\
& \leq &  \max_{y: S(y) \in    n(p \pm \delta)} || P^{t}(y,\cdot) - \pi || +
          P_x \{S(X_{\alpha n}) \notin n(p \pm \delta)\} .
\end{eqnarray*}
The last term on the right hand side doesn't depend on $t$ and is in fact in $o(1)$, uniformly in $x \in \mathcal{X}$, because when writing $S_{\alpha n} := S(X_{\alpha n})$ we get for large $\alpha$ 
\begin{eqnarray*}
P_x \{S(X_{\alpha n}) \notin n(p \pm \delta)\} & = & P_x \left\{ | S_{\alpha n} - E_x S_{\alpha n} + E_x S_{\alpha n} - n p \, | \geq \delta n \right\} \\
& \leq & P_x \left\{ | S_{\alpha n} - E_x S_{\alpha n} | \geq \delta n - | E_x S_{\alpha n} - np \, | \right\} \\
& \leq & P_x \left\{ | S_{\alpha n} - E_x S_{\alpha n} | \geq n \left( \delta - \exp{\left\{-\frac{1+\theta}{2} \alpha\right\}} \right) \right\} \\
& \leq & 2 \exp{\left\{- \frac{2}{9} \frac{n^2 \left( \delta - \exp{\left\{-\frac{1+\theta}{2} \alpha \right\}} \right)^{\!2}}{\OO{n}} \right\}} \\
& = & o(1) .
\end{eqnarray*}
Here the second inequality holds since we get from Proposition \ref{exp_theta} that 
$|E_x S_{\alpha n} - np| \leq n \gamma^{\alpha n} \leq n \exp{\{-\frac{1+\theta}{2} \alpha\}}$. 
The third inequality 
follows from our concentration result based on Azuma's inequality (Proposition \ref{conc_theta}). This proves (\ref{burnin_1}).
So after an initial $\alpha n$ steps, we may assume that the number of ones is within $\delta n$ of its stationary mean $np$.

\subsection{Coupling}
Now we project down to our two-dimensional chain: Fix any $y \in \{0,1\}^n$ such that $k:=S(y) \in n(p \pm \delta)$. Writing the transition kernel in the parametrization (\ref{transition_kernel_2dr}), we get from Corollary \ref{proj_2d} that
\begin{eqnarray}
\label{burnin_1a}
|| P^{t}(y,\cdot) - \pi || & = & || \mathcal{D}(_{S(y)}Z_{t}) - \pi_{Z_y} || \nonumber \\
& = & || P^{t}((-k,k),\cdot) - \pi_{Z_y} || \nonumber \\
& \leq & \max_{(v,v')} || P^{t}\big((v,v'),\cdot\big) - \pi_{Z_y} || \\
& \leq & \max_{(v,v'),(w,w')} || P^{t}\big((v,v'),\cdot\big) - P^{t}\big((w,w'),\cdot\big) || \nonumber \\
& \leq & \max_{(v,v'),(w,w')}    P_{(v,v'),(w,w')} \left\{ \tau > t \right\}  \nonumber .
\end{eqnarray}
The maxima here are over the entire state space of the two-dimensional chain $(v,v'),(w,w') \in \{-k,...,k\} \times \{k,...,2n-k\}$. As is well known, the second inequality above follows from stationarity of $\pi_{Z_y}$; see for instance \cite[Lemma 4.11 on page 53]{LevinPeresWilmer2009}. The last inequality is the coupling inequality, where $\tau := \min \{j \geq 0 \, : \, Z_j = Y_j\}$ is the coupling time in the coupling $(Z_j,Y_j)$ that we are about to describe.

\vspace{0.5em}

Fix any $(v,v'),(w,w') \in \{-k,...,k\} \times \{k,...,2n-k\}$ and set $Z_0 := (v,v')$ and $Y_0 := (w,w')$. Let $t=s+u$, where $s := \frac{1}{1+\theta} n \log n$ and $u := \alpha n$. 
Throughout all steps  $j=1,2,...,s,...,t$ we use an (alternating) independence coupling. By this we mean that at each step we flip a fair coin
 to decide which chain to move according to the non-lazy version of its transition kernel. The other chain stays at its current location.
Here, if $P$ is the transition probability matrix (\ref{transition_kernel_2dr}) of the chain $(Z_j)$, then $P' := 2P-I$ is its non-lazy version, where $I$ is the identity. 
However, if $Z_j$ and $Y_j$ ever agree in the $r$ (or $r'$) coordinate, we modify the coupling so that they agree in that coordinate forever after. This is possible since the probability of moving up (or down) in the $r$-coordinate does not depend on the current location in the $r'$-coordinate. Similarly, the probability of moving up (or down) in the $r'$-coordinate does not depend on the current location in the $r$-coordinate. This is easily seen from the transition kernel (\ref{transition_kernel_2dr}). 

We could implement this change  as follows: Suppose $Y_j^{(r)}=Z_j^{(r)}$. Flip fair coin number one; if it comes up heads, try to move $Y_j$ according to its non-lazy transition rule. If that would result in $Y_j$ moving up (or down) in the $r$-coordinate, flip another fair coin. If it comes up heads, move $Y_j^{(r)}$ accordingly and move $Z_j^{(r)}$ in the same way; if it comes up tails, reject the move. If fair coin number one comes up tails, try to move $Z_j$ according to its non-lazy transition rule. If that would result in $Z_j$ moving up (or down) in the $r$-coordinate, flip another fair coin. If it comes up heads, move $Z_j^{(r)}$ accordingly and move $Y_j^{(r)}$ in the same way; if it comes up tails, reject the move. Similarly for $Y_j^{(r')}=Z_j^{(r')}$.

\vspace{0.5em}

For steps $j=1,...,s$ we exploit the drift of both chains towards their stationary mean and use our calculation of the expected location of the two chains after $s$ steps. From Remark \ref{end_stage_one} we know that at time $s$ the two chains will be within $\OO{\sqrt{n}\,}$ of their stationary means in both coordinates in expectation.
For steps $j=s+1,...,t$, the drift is likely to be weak, so we can compare our two chains to two simple random walks on $\mathbb{Z}^2$. We want to show that in additional $\alpha n$ steps the two chains will coalesce with high probability.
In \cite{Barta2012} we construct an explicit coupling of $(Z_j,Y_j)$ with two lazy simple random walks $(V_j,W_j)$ to achieve this. Here we use a more abstract result due to Levin, Peres and Wilmer \cite[Proposition 17.20 on page 240]{LevinPeresWilmer2009} that is based on martingale theory.
\begin{prop}[{\bfseries Levin, Peres, Wilmer, 2009}]
\label{superMG}
Let $(R_j)_{j \geq 0}$ be a non-negative supermartingale, adapted to the sequence $(W_j)$, and let $\tau$ be a stopping time for the sequence $(W_j)$. 
Suppose that
\begin{enumerate}[(i)]
\item $R_0 = k$ ,
\item $R_{j+1} - R_j \leq B$ ,
\item there exists a constant $\sigma^2 > 0$ such that ${\rm Var}(R_{j+1} \, | \, W_0,W_1,...,W_j) \geq \sigma^2$ on the event $\{\tau > j\}$.
\end{enumerate}
If $u > 12B^2/\sigma^2$, then
$$
P_k \{\tau > u\} \leq \frac{4 k}{\sigma \sqrt{u}} .
$$
\end{prop}
We will apply this result to the two difference processes formed by
$$
R_j := Z_j^{(r)} - Y_j^{(r)} \;\;\; \mbox{ and } \;\;\; R_j' := Z_j^{(r')} - Y_j^{(r')} .
$$
It will later be convenient to restart time at $s$, so let
\begin{eqnarray*}
\tau_r := \inf \{j \geq 0 : Z_{s+j}^{(r)} = Y_{s+j}^{(r)} \} , \\
\tau_{r'} := \inf \{j \geq 0 : Z_{s+j}^{(r')} = Y_{s+j}^{(r')} \} , 
\end{eqnarray*}
so that on the event $\{\tau>s\}$, the coupling time is
\begin{equation}
\label{tau_tau_r_tau_r_prime}
\tau := \inf \{i \geq 0 : Z_i = Y_i \} = s + (\tau_r \vee \tau_{r'}) \leq s + \tau_r + \tau_{r'} .
\end{equation}
To show that the assumptions of Proposition \ref{superMG} are satisfied by $(R_j)$ and $\tau_r$, 
%(respectively $(R_j')$ and $\tau_{r'}$) 
starting at time $j=s$, we use the sequence $W_j := (Z_j,Y_j)$ and write $\mathcal{F}_j$ for the sigma algebra generated by $(W_i : i=0,1,...,j)$.
First note that, without loss of generality, we may assume that $R_s \geq 0$. If not, simply replace $R_j$ by $-R_j$. By construction of our coupling, it then follows that $R_j \geq 0$ for all $j \geq 0$.
Furthermore, the increments of $(R_j)$ are bounded above by $B := 2$,
\begin{eqnarray*}
R_{j+1} - R_j & = & \left( Z_{j+1}^{(r)} - Z_j^{(r)} \right) - \left( Y_{j+1}^{(r)} - Y_j^{(r)} \right) \\
& \leq & 2 ,
\end{eqnarray*}
since in our coupling only one of the two chains $(Z_j)$ and $(Y_j)$  can move at each step, and their (absolute) step sizes are  both bounded by two.

To see that $(R_j)$ is a supermartingale, let 
$$
D(Z_j^{(r)}) := E [Z_{j+1}^{(r)} - Z_j^{(r)} \; | \; Z_j]
$$
 be the drift of the coordinate $Z_j^{(r)}$ at time $j$, and similarly for $Y_j^{(r)}$. 
Then 
$$
E[R_{j+1} - R_j \; | \; \mathcal{F}_j]  =  D(Z_j^{(r)}) - D(Y_j^{(r)})  \leq  0 .
$$
This follows, since by (\ref{drift_function}), the drift function $D(\cdot)$ is  decreasing (non-increasing), 
and  we have $Z_j^{(r)} \geq Y_j^{(r)}$ by assumption. 
 It remains to show that the conditional one-step variances of $(R_j)$ are bounded away from zero.
Since $A := R_{j+1} - R_j \in \{-2,0,2\}$, we have to lower bound
\begin{eqnarray}
\label{var_bound_R}
& & \hspace{-3em} \frac{1}{4} {\rm Var} [R_{j+1} \; | \; \mathcal{F}_j] \nonumber \\ 
& = & \frac{1}{4} {\rm Var} [R_{j+1} - R_j \; | \; \mathcal{F}_j] \\
& = & \left[ \tilde{P}(A=2) + \tilde{P}(A=-2) \right]  -  \left[ \tilde{P}(A=2) - \tilde{P}(A=-2) \right]^2 \nonumber \\
& = & \tilde{P}(A=2) \tilde{P}(A \neq 2) + \tilde{P}(A=-2) \tilde{P}(A \neq -2) + 2 \tilde{P}(A=2) \tilde{P}(A=-2) . \nonumber
\end{eqnarray}
Here we write $\tilde{P} := P(\cdot \; | \; \mathcal{F}_j)$. From the transition probabilities (\ref{transition_kernel_2dr}) we get
$$
\tilde{P}(A=2) = \frac{1}{4n} \left[\left( k-Z_j^{(r)} \right) + \left( k + Y_j^{(r)} \right) \theta \right] ,
$$
and 
$$
\tilde{P}(A=-2) = \frac{1}{4n} \left[\left( k+Z_j^{(r)} \right) \theta + \left( k - Y_j^{(r)} \right)  \right] .
$$
From the burn-in (\ref{burnin_1}), we may assume $k \in n(p \pm \delta)$, where $\delta > 0$ is a small constant. Since $Z_j^{(r)}, Y_j^{(r)} \in [-k,k]$, this implies that $\tilde{P}(A=2)$ and $\tilde{P}(A=-2)$ are both bounded above by
$\frac{k}{n} \frac{1+\theta}{2} \leq p+\delta$,
and therefore $\tilde{P}(A \neq 2) \geq 1 - p - \delta$ and $\tilde{P}(A \neq -2) \geq 1 - p - \delta$. Consequently, the 
%quantity
%$\tilde{P}(A=2) \tilde{P}(A \neq 2) + \tilde{P}(A=-2) \tilde{P}(A \neq -2) + 2 \tilde{P}(A=2) \tilde{P}(A=-2)$ 
right hand side in (\ref{var_bound_R})
is bounded below by
$$
\left( 1-p-\delta \right) \, \left( \tilde{P}(A=2) + \tilde{P}(A=-2) \right) .
$$
Since we have the bound
\begin{eqnarray*}
\tilde{P}(A=2) + \tilde{P}(A=-2) & = & 
\frac{1}{4n} \left[ 2k(1+\theta) - (1-\theta) \left( Z_j^{(r)} + Y_j^{(r)} \right) \right] \\
& \geq & \frac{1}{4n} \left[ 2k(1+\theta) + 2k(1-\theta) \right] \\
& \geq & \frac{k}{n} \theta \\
& \geq & (p - \delta) \theta ,
\end{eqnarray*}
we get 
$$
{\rm Var}(R_{j+1} \; | \; \mathcal{F}_j) \geq \sigma^2 > 0 ,
$$
for $\sigma^2 := 4 (1-p-\delta) (p-\delta) \theta$, as required.

Analogously it can be shown that $(R_j')$ and $\tau_{r'}$ also satisfies the assumptions of Proposition \ref{superMG} 
with the same $\sigma^2$. 
Since a modification of this argument will be needed  in section \ref{sec_theta_n}, we briefly sketch it here.
Write $A' := R_{j+1}' - R_j'$, so we have to bound
$$
\tilde{P}(A'=2) = \frac{1}{4n} \left[ \left( 2n - (Z_j^{r'}+k) \right) \theta + Y_j^{r'} - k  \right]
$$
and 
$$
\tilde{P}(A'=2) = \frac{1}{4n} \left[ \left( 2n - (Y_j^{r'}+k) \right) \theta + Z_j^{r'} - k  \right] .
$$
Since  $Z_j^{r'}, Y_j^{r'}$  always  lie in $[k,2n-k]$, we
get that $\tilde{P}(A'=2)$ and $\tilde{P}(A'=-2)$ are both bounded above by $\frac{n-k}{n} \frac{1+\theta}{2} \leq 1-(p-\delta)$.
For the lower bound, we get
\begin{eqnarray*}
\tilde{P}(A'=2) + \tilde{P}(A'=-2) 
& = &  \frac{1}{4n} \left[ (Z_j^{r'} + Y_j^{r'})(1-\theta) - 2k(1+\theta) + 4n \theta \right] \\
& \geq & \frac{1}{4n} \left[ 2k(1-\theta) - 2k(1+\theta) + 4n \theta \right] \\
& = & \frac{n-k}{n} \theta \\
& \geq & (1 - p -\delta) \theta .
\end{eqnarray*}
In the same way as above for $(R_j)$, this shows that ${\rm Var}(R_{j+1}' | \mathcal{F}_j) \geq \sigma^2 := 4(1-p-\delta)(p-\delta)\theta > 0$ in this case as well.

Applying Proposition \ref{superMG}, we get 
%that $R_j := Z_j^{(r)} - Y_j^{(r)}$ and $R_j' :=  Z_j^{(r')} - Y_j^{(r')}$ satisfy (for $u > 48/\sigma^2$)
for $u > 48/\sigma^2$, 
\begin{eqnarray}
\label{entropy_result}
P_{Z_s,Y_s} \{ \tau_r > u \} & \leq & \frac{4 |R_s|}{\sigma \sqrt{u}}  \;\;\; \mbox{ and} \\
P_{Z_s,Y_s} \{ \tau_{r'} > u \} & \leq & \frac{4 |R_s'|}{\sigma \sqrt{u}} . \nonumber
\end{eqnarray}

\vspace{0.5em}

Now we put things together. By conditioning on where we are after the first $s$ steps (the drift regime of our coupling), we get
\begin{eqnarray*}
P_{(v,v'),(w,w')} \{ \tau > s + u \, | \, Z_s, Y_s\} & = & \mathbf{1}\{\tau > s\} P_{(v,v'),(w,w')} \{ \tau > u \; | \; Z_s,Y_s\} \\
& \leq & P_{Z_s,Y_s} \{ \tau_r + \tau_{r'} > u \} \\
& \leq & P_{Z_s,Y_s} \{\tau_r > u/2\} + P_{Z_s,Y_s} \{\tau_{r'} > u/2\} \\
%& \leq & P_{R_s} \{\tau_r > u/2\} + P_{R_s'} \{\tau_{r'} > u/2\} \\
%& \leq & P_{D_s^{(r)}} \{\tilde{\tau}_r > u/2\} + P_{D_s^{(r')}} \{\tilde{\tau}_{r'} > u/2\} \\
%& \leq & \frac{C | D_s^{(r)} |}{\sqrt{\delta_r u}} + \frac{C | D_s^{(r')}|}{\sqrt{\delta_{r'} u}}  .\\
& \leq & \frac{4 \sqrt{2}}{\sigma} \;  \frac{|R_s| + |R_s'|}{\sqrt{u}} .
\end{eqnarray*}
Here the first inequality comes from (\ref{tau_tau_r_tau_r_prime}), where we restart time at $s$.
The last inequality is (\ref{entropy_result}).
By taking expectation, we get for some positive constants $C,C'$ and large $n$ that
\begin{eqnarray*}
P_{(v,v'),(w,w')} \{ \tau > s+u \} & \leq & 
\frac{4 \sqrt{2}}{\sigma \sqrt{u}} \; \left[ E_{(v,v'),(w,w')} |R_s| + E_{(v,v'),(w,w')} |R_s'| \right] \\
& = & \frac{4 \sqrt{2}}{\sigma \sqrt{u}} \; \left[ |E_{(v,v'),(w,w')} R_s| + |E_{(v,v'),(w,w')} R_s'| \right] \\
& \leq & \frac{4 \sqrt{2}}{\sigma} \; \frac{C \sqrt{n} + C' \sqrt{n}}{\sqrt{\alpha n}} .
\end{eqnarray*}
For the equality above, note that we can move the absolute values outside the expectations, since in our coupling the sign of $R_j$ (respectively $R_j'$) can never switch between plus and minus.
The second inequality above comes from Remark \ref{end_stage_one}.
Combining this with (\ref{burnin_1}) and (\ref{burnin_1a}), we get 
\begin{eqnarray*}
d(t) & := & \max_x ||P^t(x,\cdot) - \pi|| \\
& \leq & \max_{(v,v'),(w,w')} P_{(v,v'),(w,w')} \{\tau > t\} + \oo{1} \\
& \leq & \frac{4 \sqrt{2} (C+C')}{\sigma \sqrt{\alpha}} + \oo{1} ,
\end{eqnarray*}
which implies
$$
\lim_{\alpha \rightarrow \infty} \limsup_{n \rightarrow \infty} d(t) = 0 .
$$
This finishes the proof of the upper bound part of Theorem \ref{theta_model}. $\Box$

\section{Metropolis vs. Gibbs }
\label{RWM_GS}
Since the coordinates are independent under the distribution $\pi(x) = \theta^{S(x)} (1+\theta)^{-n}$ on $\mathcal{X}=\{0,1\}^n$, where $\theta \in (0,1]$, 
there exists a very close connection between the random walk Metropolis algorithm for $\pi$, that we studied above, and the Gibbs sampler for the distribution $\pi$.
The (random scan) Gibbs sampler for the distribution $\pi$, also known as Glauber dynamics in statistical physics, is also a Markov chain $(X_t)$ on $\mathcal{X}$ that converges to $\pi$ as the number of steps $t$ goes to infinity. It evolves as follows: Given we are at $X_t=x$, draw a coordinate $i \in [n]$ uniformly at random. Then set $X_{t+1}^{(j)} := X_t^{(j)}$ for all $j \neq i$, and set $X_{t+1}^{(i)} := B$, where $B \in \{0,1\}$ is a draw from the conditional distribution of $\pi$ for coordinate $i$, given the values of $x$ at all other coordinates $j \neq i$. 
%I.e.
%$$
%B \sim \pi(X^{(i)} \, | \, X^{(j)}=x^{(j)} \mbox{ for all } j \neq i) .
%$$
Since the coordinates are all i.i.d. Ber$(p)$ under $\pi$, where $p=\frac{\theta}{1+\theta}$, the random variable $B$ above  follows a Ber$(p)$ distribution and is independent of 
$X_t$.
%the values $x^{(j)}$ at the other coordinates $j \neq i$.
\begin{defin}
Let $\mathbb{P}$ be the transition probability matrix of the (non-lazy) random walk Metropolis algorithm for $\pi$.
Similarly, let $\mathbb{Q}$ be the transition probability matrix of the (non-lazy) Gibbs sampler for $\pi$.
Fix any $q \in [0,1]$. Then the \emph{lazy(q) random walk Metropolis algorithm for $\pi$} refers to the Markov chain on $\mathcal{X}$ with transition probability matrix $(1-q) \mathbb{P} + q I$. Similarly, the \emph{lazy($q$) Gibbs sampler for $\pi$} refers to the Markov chain on $\mathcal{X}$ with transition probability matrix $(1-q) \mathbb{Q} + q I$. Here $I$ is the (appropriately sized) identity matrix. 
We refer to $q$ as the laziness-factor of the chain in question.
When no value $q$ is specified, 
a lazy chain is understood to be lazy$(1/2)$.
\end{defin}

\begin{remark}
The lazy$(q)$ random walk Metropolis algorithm for $\pi$ simply refers to the Markov chain on $\mathcal{X}$, where at each step with probability $q$ we stay where we are, and with probability $1-q$ we make one step according to the (non-lazy) random walk Metropolis algorithm. Similarly for the Gibbs sampler. 
For example, the lazy random walk Metropolis algorithm studied before corresponds to the lazy($1/2$) random walk Metropolis algorithm, whereas the lazy($0$) Gibbs sampler refers to the non-lazy Gibbs sampler for $\pi$.
\end{remark}
With this definition, we get the following relationship between the random walk Metropolis algorithm and the Gibbs sampler for $\pi$.

\begin{prop}
\label{Metropolis_Gibbs}
Fix any $\theta \in (0,1]$ and consider the distribution $\pi=\pi_{\theta}$ on $\mathcal{X}=\{0,1\}^n$ defined by $\pi(x) = \theta^{S(x)} (1+\theta)^{-n}$.
Then the lazy $\!(1/2)$ random walk Metropolis algorithm for $\pi$ follows the same law as the lazy $\!\!\big((1-\theta)/2\big)$ Gibbs sampler for $\pi$.
\end{prop}

\begin{remark}
For $\theta = 1$ this shows that the lazy random walk Metropolis algorithm for the uniform distribution on the hypercube is equal to the non-lazy Gibbs sampler for the uniform distribution. 
If we say that a non-lazy chain ``runs twice as fast'' as its lazy ($q=1/2$) version, then for $\theta=1$
the lazy random walk Metropolis algorithm ``runs twice as fast'' as the lazy Gibbs sampler. 
The above Proposition shows that this ``speed advantage'' of the Metropolis algorithm disappears as $\theta \downarrow 0$.
The reason for this is that the random walk proposal in the Metropolis algorithm is less and less well adapted to $\pi$ as $\theta$ goes down to zero.
\end{remark}

\noindent  \textbf{Proof of Proposition \ref{Metropolis_Gibbs}.} Fix any $\theta \in (0,1]$. The lazy($\frac{1}{2}$) random walk Metropolis algorithm for $\pi$ has the following transition kernel:
\begin{equation} 
\label{transition_kernel_lRW-MH-theta}
\mathbb{P}(x,y) =  \left\{ \begin{array}{l@{\quad:\quad}l}
\frac{1}{2n} & x \sim y, S(x) > S(y)  , \\
\frac{1}{2n} \theta & x \sim y, S(x) < S(y) , \\
\frac{1}{2} + \frac{n-S(x)}{2n} (1-\theta) & x = y , \\
0 & \mbox{otherwise}  .\\ 
\end{array} \right.
\end{equation}
%
%
%We need to show that this corresponds exactly to the transition kernel $\mathbb{Q}(x,y)$ for the lazy($\frac{1-\theta}{2}$) Gibbs sampler for $\pi$.
%
%
As noted before, it's easy to see that we can run the corresponding  chain $(X_t)$ in the following way. Given we are at state $X_t$ at time $t$:
\begin{itemize}
\item Pick a coordinate $i \in [n]$ uniformly at random, independent of all previous choices.
\item Draw $U_t \sim$ Uniform$[0,1]$, independent of all previous choices.
\item Set $X_{t+1}^{(j)} := X_t^{(j)}$ for $j \ne i$, and set the $i^{th}$ coordinate of $X_{t+1}$ to
\end{itemize}
$$
X_{t+1}^{(i)} :=  \left\{ \begin{array}{l@{\quad:\quad}l}
1 & 0 \leq U_t \leq \frac{\theta}{2} , \\
X_t^{(i)} & \frac{\theta}{2} < U_t \leq \frac{1}{2} , \\
0 & \frac{1}{2} < U_t \leq 1   .\\
\end{array} \right.
$$
However, this corresponds exactly to the transition rule of the lazy$(\frac{1-\theta}{2})$ Gibbs sampler: If $U_t \in (\frac{\theta}{2},\frac{1}{2}]$, that is, with probability $(1-\theta)/2$, we stay where we are; and if $U_t \notin (\frac{\theta}{2},\frac{1}{2}]$, that is, with probability $1 - (1-\theta)/2 = (1+\theta)/2$, we pick one of the $n$ coordinates uniformly at random and replace it with an independent Bernoulli$(\frac{\theta}{1+\theta})$ random variable, since we get $P(U_t \leq \theta/2 \, | \, U_t \notin (\frac{\theta}{2},\frac{1}{2}]) = \frac{\theta}{1+\theta}$.
$\Box$
%
%
%
\begin{comment}
%
% Alternative proof:
%
So fix any $x \in \mathcal{X}$.
In the Gibbs sampler, the only possible transitions from $X_t=x$ are to one of the $n$ neighboring states $y \sim x$, or to stay at $x$. 
Pick any $y \sim x$ such that $x^{(i)}=0, y^{(i)}=1$ for some $i \in [n]$, and (necessarily) $x^{(j)}=y^{(j)}$ for all $j \neq i$. Then the transition probability from $x$ to $y$ in the lazy($\frac{1-\theta}{2}$) Gibbs sampler for $\pi$ is
\begin{eqnarray*}
\mathbb{Q}(x,y) & = & P(\mbox{make a non-lazy move}) \, P(\mbox{pick coordinate } i \, | \, \mbox{move}) P(B = 1 \, | \, \mbox{move}) \\
& = & \left(1 - \frac{1-\theta}{2}\right) \, \frac{1}{n} \, \frac{\theta}{1+\theta} \\
& = & \frac{1}{2n} \theta \\
& = & \mathbb{P}(x,y) .
\end{eqnarray*}
Similarly, pick any $y \sim x$ such that $x^{(i)}=1, y^{(i)}=0$ for some $i \in [n]$, and (necessarily) $x^{(j)}=y^{(j)}$ for all $j \neq i$.
Then we get
$$
\mathbb{Q}(x,y) = \left(1 - \frac{1-\theta}{2}\right) \frac{1}{n} \left(1 - \frac{\theta}{1+\theta}\right) = \frac{1}{2n} = \mathbb{P}(x,y) .
$$
So the two transition probability matrices are in fact the same. $\Box$
%
\end{comment}
%
%
%

\section{Edge probabilities varying with $n$}
\label{sec_theta_n}
So far we considered the model $\pi(x) = \theta^{S(x)} (1+\theta)^{-n}$ where $\theta \in (0,1]$ is fixed as $n$ goes to infinity. In this section we consider the case where $\theta = \theta_n$ is allowed to vary with $n$. As mentioned in \cite[page 2118]{DiaconisSaloff-Coste2006}, 
the methods and results of \cite{DiaconisRam2000} imply that
the non-lazy version of the Metropolis algorithm for this model has cutoff at $(1+\theta_n)^{-1} n \min\{ \log n, \log \sqrt{n / \theta_n}\}$ with window size $n$. 
By the same argument as for Theorem \ref{theta_model}, this implies the following result for the lazy version:
\begin{theorem}[{\bfseries Diaconis, Ram, 2000}]
\label{theta_n_model}
The lazy random walk Metropolis chain for $\pi(x) = \theta_n^{S(x)} (1+\theta_n)^{-n}$ on $\{0,1\}^n$, where $\theta_n \in (0,1]$, has cutoff at $\frac{2}{1+\theta} n \min\{ \log n, \log \sqrt{n / \theta_n}\}$
 with a window of size $n$.
\end{theorem}
\begin{corollary}
Let $n := {\nu \choose 2}$ and let $\pi(x) = p_n^{S(x)} (1-p_n)^{n-S(x)}$ be the Erd\H{o}s-R\'enyi random graph model on
$\nu$ vertices with parameter $p_n \in (0,1)$. 
The lazy random walk Metropolis chain for this model  has cutoff at 
$$
2 \mspace{0.5mu} \rho \mspace{1.5mu} n \min\left\{\log n, \log \sqrt{\frac{n \rho}{1-\rho}}\right\}
$$
 with a window of size $n$, where $\rho := \max \{p, 1-p\}$.
\end{corollary}
\begin{remark}
By absorbing $- (1+\theta)^{-1} n \log \theta$ into the window size term, Theorem \ref{theta_n_model} is seen to be consistent with the result from Theorem \ref{theta_model} for constant $\theta$. 
We note that the expression for the mixing time changes when $\theta_n \leq 1/n$, corresponding to $p_n := \theta_n/(1+\theta_n) \leq 1/(n+1)$,
which is right about where the edge probability $p_n$ passes the critical threshold of $1/n$ for the existence of a giant component in the Erd\H{o}s-R\'enyi random graph.
\end{remark}

\noindent  \textbf{Proof of Theorem \ref{theta_n_model}.} 
In the remainder of this section we will extend our proof of Theorem \ref{theta_model} to the case where $\theta$ is allowed to vary with $n$, 
and prove Theorem \ref{theta_n_model}. 
For the lower bound part 
of this cutoff result, 
we can follow the proof given in section \ref{sec_lower_bound_1}. Propositions \ref{exp_theta} and \ref{coupon_collector} including their proofs stay true as stated when $\theta_n$ is allowed to vary with $n$.
After that, use the bound
$$
{\rm Var}_{\underline{1}} S_t \leq \frac{n}{(1+\theta_n)^2} \left[ \gamma^t + \theta_n \right] ,
$$
which leads to the bound
$$
\sigma^2 := 1/2 \left( {\rm Var}_{\underline{1}} S_t + {\rm Var}_{\pi_S} S \right) \leq 
\frac{n}{(1+\theta_n)^2} \left[ \gamma^t + \theta_n \right] .
$$
Then we get
\begin{eqnarray}
|E_{\underline{1}} S_t - E_{\pi_S} S| & = & \frac{n}{1+\theta_n} \gamma^t  \nonumber \\
& \geq & \sigma \frac{\sqrt{n} \gamma^t}{\sqrt{\gamma^t + \theta_n}} . \label{expectation_difference_theta_n}
\end{eqnarray}
For $t := \frac{2}{1+\theta_n} n \min\{ \log n, \log \sqrt{n/\theta_n}\} - \alpha n$ and $\gamma := 1 - \frac{1+\theta_n}{2n}$ as before, we get
$$
\gamma^t \sim \max\{ 1/n, \sqrt{\theta_n/n} \} \; e^{\frac{1+\theta_n}{2} \alpha} ,
$$
which implies
$$
\lim_{\alpha \rightarrow \infty} \liminf_{n \rightarrow \infty} \left( \frac{\sqrt{n} \gamma^t}{\sqrt{\gamma^t + \theta_n}} \right)^{\!2} = \infty .
$$
%
% Proof:
%
% Case 1: $\theta_n \leq 1/n$. Note that $t=2/(1+\theta_n) n \log n - \alpha n$ and $n \gamma^t \sim e^{\alpha (1+\theta)/2}$.
% Then $r^2 := \gamma^{2t} \frac{n}{\gamma^t + \theta_n} \geq \frac{(n \gamma^t)^2}{1 + n \gamma^t} \sim e^{\alpha (1+\theta)/2} \rightarrow \infty$.
%
% Case 2: $\theta_n > 1/n$. Note that $t=2/(1+\theta_n) n \log n - \alpha n$ and $n \gamma^t \sim \sqrt{n \theta_n} e^{\alpha (1+\theta_n)/2}$.
% Then r^{-2} := \frac{n \gamma^t + n \theta_n}{(n \gamma^t)^2} \rightarrow 0$. $\Box$
%
%
The lower bound result now follows from (\ref{expectation_difference_theta_n}) and Proposition \ref{dist_stat}. 
%$\Box$

For the upper bound part of the cutoff result, we distinguish two cases, according to which term attains the minimum in the expression for  
the mixing time.
For $\theta_n \leq 1/n$, it is enough to consider a coordinate wise coupling: To update two copies $(X_t), (Y_t)$ of the chain, use the construction given before Proposition \ref{coupon_collector} to run the chains, but update \emph{the same} coordinate $i \in [n]$ and use \emph{the same} uniform random variable $U_t$ \emph{for both} $(X_t)$ and $(Y_t)$. If $R_t$ is the number of coordinates not refreshed by time $t$, as defined before Proposition \ref{coupon_collector}, then for for $t := \frac{2}{1+\theta_n} n \log n + \alpha n$ we get
\begin{eqnarray*}
d(t) & = & \max_{x} || P^t(x,\cdot) - \pi || \\
& \leq & \max_{x,y} || P^t(x,\cdot) - P^t(y,\cdot) || \\
& \leq & \max_{x,y} P_{x,y}(X_t \neq Y_t) \\
& \leq &  P (R_t \geq 1) \\
& \leq & E R_t \\
%& \leq & \max_{x,y} P_{x,y} (R_t \geq 1) \\
%& \leq & E_{x,y} R_t \\
& = & n \gamma^t \\
& \sim & e^{-\frac{1+\theta_n}{2} \alpha} , 
\end{eqnarray*}
which proves that $\lim_{\alpha \rightarrow \infty} \limsup_{n \rightarrow \infty} d(t) = 0$, as required.

For $\theta_n > 1/n$, this coordinate wise coupling is generally not sharp enough to establish cutoff. 
To see this, note that for $t=\frac{2}{1+\theta_n} n \log \sqrt{n/\theta_n} + \alpha n$ we get $n \gamma^t \sim \sqrt{n \theta_n} e^{-\alpha (1+\theta_n)/2}$, which grows without bound unless $\theta_n \in \OO{1/n}$.
However, a slight modification of our proof for constant $\theta$ works in this case. 
With the exception of Remark \ref{end_stage_one}, 
everything in this proof stays true until the application of Proposition \ref{superMG}.
Restarting time after a burn-in of $\alpha n$ steps, we may assume that the number of ones of the state $X_0$ that we use for the two dimensional projection satisfies $k := S(X_0) \in n (p_n \pm \delta)$, for some small constant $\delta > 0$.
Again, we write $p_n := \theta_n/(1+\theta_n)$.
Note that $p_n-\delta$ can be negative in this case.
The two difference processes formed by
$$
R_j := Z_j^{(r)} - Y_j^{(r)} \;\;\; \mbox{ and } \;\;\; R_j' := Z_j^{(r')} - Y_j^{(r')} 
$$
still satisfy all the assumptions of Proposition \ref{superMG}, except that the conditional one step variances can no longer be bounded away from zero. Therefore, here we modify the argument given for constant $\theta$, treating the $r$ and $r'\!$-coordinates separately.

The width of the domain $[k,2n-k]$ of the $r'\!$-coordinate of our chain $(Z_t)=\left(Z_t^{(r)}, Z_t^{(r')}\right)$ is of order $n$, as before, and we can still apply Proposition \ref{superMG}.
To bound the conditional one step variances, using the same notation as above,  we get that $\tilde{P}(A'=2)$ and $\tilde{P}(A'=-2)$ are both bounded above by $\frac{n-k}{n} \frac{1+\theta}{2} \leq 1-\delta'$, for some constant $\delta'>0$ that doesn't depend on $n$. To see this inequality, first note that both terms on the left hand side are in $[0,1]$.
%and since $k \in n(p \pm \delta)$, if one of the terms is large, the other one is small. 
Let $\delta'' := 2\delta/(1-\delta)$, where $\delta>0$ is the constant we used in the burn-in to ensure that the number of ones $k$ of the state used in the two dimensional projection satisfies $k \in n(p \pm \delta)$.
%
% Note that $p-\delta$ might be negative here, since $p=p_n=\theta_n/(1+\theta_n)$ might go to zero as $n$ goes to infinity.
%
If $\theta_n < \delta''$, we get $\frac{1+\theta_n}{2} < \frac{1+\delta''}{2}$, so we can pick $\delta' := (1-\delta'')/2$.
If $\theta_n \geq \delta''$,
 we get $p_n - \delta \geq \frac{\delta''}{1+\delta''} - \delta > 0$, so here we can pick
$\delta' := \frac{\delta''}{1+\delta''} - \delta > 0$.
For the lower bound, we get as before
%\begin{eqnarray*}
%\tilde{P}(A'=2) + \tilde{P}(A'=-2) 
%& = &  \frac{1}{4n} \left[ (Z_j^{r'} + Y_j^{r'})(1-\theta_n) - 2k(1+\theta_n) + 4n \theta_n \right] \\
%& \geq & \frac{1}{4n} \left[ 2k(1-\theta_n) - 2k(1+\theta_n) + 4n \theta_n \right] \\
%& = & \frac{n-k}{n} \theta_n \\
%& \geq & (1 - p - \delta) \theta_n .
%\end{eqnarray*}
$$
\tilde{P}(A'=2) + \tilde{P}(A'=-2) \geq (1 - p - \delta) \theta_n .
$$
That means our variance bound ${\rm Var}(R_{j+1}' | \mathcal{F}_j) \geq \sigma^2_n := 4 \delta' (1-p-\delta) \theta_n$  depends on $n$ and can't be bounded away from zero, as before.
However, for some constant $C>0$ and large $\alpha$ we still get $u := \alpha n > C/\theta_n = 12 B^2 / \sigma_n^2$, as required, since we have $\theta_n > 1/n$. Applying Proposition \ref{superMG} to $(R_j')$ and $\tau_{r'}$ entails
\begin{equation}
\label{r'_mod1}
P_{Z_s,Y_s} \{ \tau_{r'} > u/2 \} \leq \frac{4 |R_s'|}{\sigma_n \sqrt{u/2}} .
\end{equation}
Furthermore, for $s := 2/(1+\theta_n) n \log \sqrt{n/\theta_n}$ we get $n \gamma^s \sim \sqrt{n \theta_n}$, and therefore, for some constant $c>0$,
\begin{equation}
\label{r'_mod2}
\frac{|E R_s'|}{\sigma_n \sqrt{u}} \leq \frac{n \gamma^s}{\sqrt{c \theta_n \alpha n}} \sim \frac{1}{\sqrt{\alpha c}} ,
\end{equation}
as required 
%for the application of Proposition \ref{superMG} 
in this context.

On the other hand, the width of the domain $[-k,k]$ for the $r$-coordinate of our chain $(Z_t)$ is not necessarily of order $n$, as it was for fixed $\theta$. 
As before, for $(R_j)$ and $\tau_r$ we get that $\tilde{P}(A=2)$ and $\tilde{P}(A+-2)$ are both 
upper bounded by $\frac{k}{n} \frac{1+\theta_n}{2} \leq p + \delta$, and further
$$
\tilde{P}(A=2) + \tilde{P}(A=-2) \geq \frac{k}{n} \theta_n .
$$
Therefore, we get the variance bound
$$
{\rm Var} \left[ R_{j+1} | \mathcal{F}_j \right] \geq \sigma_n^2 := 4 (1-p_n-\delta) \frac{k}{n} \theta_n .
$$
This violates the assumption $\alpha n =: u > 12 B^2/\sigma_n^2$ of Proposition \ref{superMG}.
Therefore,  we use here a slightly modified version of Proposition \ref{superMG}, where we replace the last sentence in its statement by the following one: If $h \geq 2B$ then for any $u>0$ we get
$$
P_k\{\tau > u\} \leq \frac{k}{h} + \frac{3kh}{u \sigma^2} . 
$$
This version directly follows from the proof or the original version as given by Levin, Peres and Wilmer \cite[page 346]{LevinPeresWilmer2009}. Minimizing the right hand side above by choosing $h:=\sqrt{u \sigma^2/3}$ leads to the original version of Proposition \ref{superMG}, which we can't apply in our context. Instead, here we choose $h:=\sqrt{\alpha n \theta_n}$.
Since $\theta_n>1/n$, for large $\alpha$ we get $h \geq 2B = 4$, as required. Applying this modified version of Proposition \ref{superMG} to $(R_j)$ and $\tau_r$, we get
\begin{eqnarray*}
P_{Z_s,Y_s} \{\tau_r > u/2\} & \leq & \frac{|R_s|}{h} + \frac{3|R_s|h}{\sigma^2_n u/2} \\
& \leq & \frac{|R_s|}{\sqrt{\alpha n \theta_n}} + \frac{6|R_s|\sqrt{\alpha n \theta_n}}{\alpha n \theta_n c k/n} \\
& = & \frac{|R_s|}{\sqrt{\alpha n \theta_n}} \left[ 1 + \frac{6n}{ck} \right] .
\end{eqnarray*}
Here, we used $u:=\alpha n$ and $\sigma_n^2 \geq c \theta_n k/n$ for some constant $c>0$. 
By setting $s:=2/(1+\theta_n) n \log \sqrt{n/\theta_n}$ and taking expectation, this leads to the bound
\begin{eqnarray*}
\frac{E_{(v,v'),(w,w')} |R_s|}{\sqrt{\alpha n \theta_n}} \left[ 1 + \frac{6n}{ck} \right] & = &
\frac{|E_{(v,v'),(w,w')} R_s|}{\sqrt{\alpha n \theta_n}} \left[ 1 + \frac{6n}{ck} \right] \\
& \leq & \frac{2k\gamma^s}{\sqrt{\alpha n \theta_n}} \left[ 1 + \frac{6n}{ck} \right] \\
& \sim & \frac{2\sqrt{n\theta_n} k/n}{\sqrt{\alpha n \theta_n}} \left[ 1 + \frac{6n}{ck} \right] \\
& \leq & \frac{2 + 12/c}{\sqrt{\alpha}} .
\end{eqnarray*}
Here, the first equality follows by the same argument as in the case for constant $\theta$. The first inequality follows from Remark \ref{exp_2dr_general}, using  the fact that the $r$-coordinates live in the domain $[-k,k]$.
The asymptotic equality comes from $n \gamma^s \sim \sqrt{n \theta_n}$.
In the same way as for constant $\theta$, we can now  combine this result for $(R_j)$ and $\tau_r$
with the corresponding results (\ref{r'_mod1}) and (\ref{r'_mod2}) for $(R_j')$ and $\tau_{r'}$
to establish the upper bound part of the cutoff result. 
This finishes the proof. $\Box$

\section*{Acknowledgments:}
I would like to thank my PhD advisor Steven Lalley for his invaluable guidance and many helpful discussions and suggestions. I also would like to thank Michael Wichura for a very careful reading of an earlier version of this paper that lead to significant improvements.

%\section*{References}
%

\end{document}